\documentclass[leqno,12pt]{article}
\usepackage[hang]{subfigure}
\usepackage{color}
\usepackage{color}
\usepackage{colortbl}
\usepackage{mystyle}
\usepackage{latexsym}
\usepackage{setspace}

\renewcommand{\epsilon}{\varepsilon}

\begin{document}

\title{Fourier analysis of serial dependence measures}

\author{ Ria Van Hecke\\
Ruhr-Universit\"at Bochum \\
Fakult\"at f\"ur Mathematik \\
44780 Bochum \\
Germany  \\
\and
Stanislav Volgushev \\
University of Toronto \\
Department of Statistical Sciences\\
Toronto, Ontario M5S 3G3 \\
Canada
\and 
Holger Dette\\
Ruhr-Universit\"at Bochum \\
Fakult\"at f\"ur Mathematik \\
44780 Bochum \\
Germany  \\
}

 \maketitle

\begin{abstract}
Classical spectral analysis is based on the discrete Fourier transform of the auto-covariances. In this paper we investigate the asymptotic properties of new frequency domain methods where the auto-covariances in the spectral density are replaced by alternative dependence measures which can be estimated by U-statistics. An interesting example is given by Kendall's $\tau$, for which the limiting variance exhibits a surprising behavior.

\end{abstract}


\section{Introduction}\label{sec1}
\def\theequation{1.\arabic{equation}}
\setcounter{equation}{0}

Over the years spectral analysis has developed into a fundamental important tool kit in  the analysis of data from a stationary time series $\{ X_t\}_{t\in \Z}$. The spectral density, defined as the discrete Fourier transform of the auto-covariances, provides a convenient way to characterize  the second order properties of a stationary sequence. Estimation of the spectral density is usually performed by smoothing  the  periodogram, that is the discrete Fourier transform of empirical auto-covariances [see for example Chapters~4 and~10 of \citet{BrockwellDavis1987}]. \\
 It is well known that this approach is not able to capture non-linear features of time series dynamics such as changes in skewness, kurtosis or dependence in the extremes. This motivated numerous authors to describe serial dependence by considering spectral densities corresponding to a family of transformations of the original time series 
[see \cite{hong1999,hong2000}, \cite{Li2008,Li2012}, \cite{Hagemann2013}, \cite{DetteEtAl2013}, \cite{birr2014}, \cite{DavisMikoschZhao2013}, \cite{DetteEtAl2016}]. Roughly speaking, these authors suggest to define a family of  spectral densities, say $\{ f(\lambda,x,y) ~|~x,y\}$, where the auto-covariances (at lag $k$) are replaced by functionals of the lag $k$-distributions $\Prob ( X_t \leq x,X_{t+k} \leq y ) $. This approach is attractive as it allows a more complete description of the serial dependence. The price for this flexibility is the calculation  of  a family of spectral densities, in contrast 
 to the classical approach, which uses only {\bf one} spectral density calculated as the discrete Fourier transform of the auto-covariances. \\
In the present paper we investigate a class of alternative spectral densities, which keeps the simplicity of the classical spectral theory but eliminates some drawbacks arising from the use of auto-covariances in its definition. More precisely, we consider general spectral densities of the form 
\als{
\spec_{\xi}(\omega)=\frac{1}{2\pi}\sum_{k \in \Z}\xi_k e^{-ik\omega}\qquad(\omega\in\R),\label{uspecdens} 
}
where  for  each $k \in \Z$ the quantity  $\xi_k$ denotes  a dependence measure between the random variables $X_t$ and $X_{t+k}$ (in the classical case $\xi_k=r_k=$Cov$(X_t,X_{t+k})$) and we implicitly assumed that $\sum_{k\in\Z}|\xi_k|<\infty$. Spectral densities  of the form \eqref{uspecdens} have been considered by \cite{ahdesmaki2005}, \cite{zhou2012} and \cite{carcea2015}, who replaced the lag $k$ auto-covariance by other measures of dependence such as Kendall's $\tau$, distance correlation, or L-moments. A thorough theoretical analysis of this idea for dependence measures that can be represented as {\bf \textit{linear}} functionals of the empirical copula at lag $k$ was conducted in \cite{DetteEtAl2016}. Their analysis includes dependence measures such as Spearman's rank autocorrelation [see \cite{WaldWolfowitz1943}], Blomqvist's beta [see \cite{blomqvist1950}] and Gini's rank association coefficient [see \citet{Schechtman1987}]. However, the theory depends crucially on the linearity of the corresponding functional and cannot be generalized to other dependence measures. A particularly interesting dependence measure that is not covered by the analysis of \cite{DetteEtAl2016} is Kendall's tau which can be represented as by  
\als{
\tau_k = 4\int C_k(u) dC_k(u) - 1
\label{kendal}
}
where $C_k$ denotes the copula corresponding to lag $k$. Note that Kendall's tau  is a non-linear functional of the lag $k$ copula. 
A classical approach to the estimation of Kendall's tau is based on the  representation
\[
\tau_k = 2\Prob[X_1<X_2,Y_1<Y_2]+2\Prob[X_2<X_1,Y_2<Y_1]-1
\]  
where $(X_1,Y_1), (X_2, Y_2)$ are independent copies with the same distribution as $(X_0,X_k)$. Motivated by this example we are interested in the statistical properties of estimators of spectral densities of the form \eqref{uspecdens} with a 
measure $\xi_k$ of lag $k$ dependence that can be represented as 
\als{\label{urepres}
\xi_k=\E\Big[h\Big(\colvec{X_{0}^{(1)}}{X_{k}^{(1)}},\dots,\colvec{X_{0}^{(m)}}{X_{k}^{(m)}}\Big)\Big].
}
where  $( {X_{0}^{(1)}},{X_{k}^{(1)}} )^T , \ldots, ( {X_{0}^{(m)}},{X_{k}^{(m)}})^T $ are 
 independent copies  of $ ({X_0}, {X_k})^T $   and  $h$ is a symmetric kernel of order $m$. The representation 
 \eqref{urepres}  motivates to estimate  $\xi_k$ by a $U$-statistic, say $\xi_{n,k}$, and to form the  corresponding 
 \textit{$U$-lag-window estimate}
\als{\label{per}
\hat{f}_{n,\xi}(\omega)=\frac{1}{2\pi}\sum_{|k|<n}w_n(k)\xi_{n,k}e^{-ik\omega},
}
where $\{w_n(k)\}_{k=-(n-1),\ldots , n-1} $ are given weights. 
In Section \ref{sec2}  we will introduce the necessary notation and illustrate the general approach by  several examples.
The main results of the paper can be found in Section \ref{sec3}, where we investigate the asymptotic properties of the new
estimates. In  particular we prove consistency of the 
estimate \eqref{per} for a broad class of kernels $h$ and establish its  asymptotic normality for several important cases including Kendalls $\tau$.
Interestingly the asymptotic variance of the $U$-lag-window estimate based on  Kendall's tau depends on the spectral density \eqref{uspecdens}
where the quantities $\xi_k$ are the lag $k$ Spearman's rho correlations.
The proofs are very involved and will be deferred to Section   \ref{sec4}, while more technical  arguments can be found in Section \ref{sec5}.

\section{Examples of U-lag-window spectral densities and their estimators} 
\label{sec2}
\def\theequation{2.\arabic{equation}}
\setcounter{equation}{0}

Throughout this paper let $\{ X_t \}_{t\in \Z}$ be a real-valued process and denote by $F$ and $F_k$ the marginal distribution function of $X_t$ and the distribution function of the pair $(X_t,X_{t+k})$, respectively ($k\in\Z$). Recall the definition of the spectral density $\spec_{\xi}$ in \eqref{uspecdens}, where the measure of dependence (at lag $k$) has the representation \eqref{urepres} for a given kernel $h$ of order $m$. Throughout this paper, we will maintain the following assumption
\bi  
\item[(C0)] The process $\{ X_t \}_{t\in \Z}$ is strictly stationary. The functions $F$ and $F_k$ are continuous (for all $k\in\Z$) and $\sum_{k\in\Z}|\xi_k|<\infty$.
\ei
Let  $X_1,\dots,X_n$ be the finite stretch of this process representing the observed data and 
define for any $k\in\{-(n-1),\dots\,n-1\}$ the set  $\mathcal{T}_k:=\big\{t|t,t+k\in\{1,\dots,n\}\big\}$. In the following example we illustrate how different kernels yield different measures 
of dependence and as consequence different spectral densities.

\begin{example}\label{exustat}{\rm \quad
\bi[leftmargin=*]
\item[(i)] If $m=2$ and    $h\big ( ({x_1},{y_2})^T , ({x_2},{y_2})^T \big) =\frac{1}{2}(x_1-x_2)(y_1-y_2)$, then 
the representation \eqref{urepres} gives the auto-covariance at lag $k$, that is 
\[
r_k = \E\Big[\frac{1}{2}(X_0^{(1)}-X_0^{(2)})(X_k^{(1)}-X_k^{(2)})\Big]=\Cov(X_0,X_k)
\]
and we obtain the classical  spectral density. 
\item[(ii)] If $m=2$, $I(\cdot)$ denotes the indicator function and the kernel is defined by 
\als{\label{kerneltau} 
h\big ( ({x_1},{y_2})^T , ({x_2},{y_2})^T \big) 
=2I(x_1<x_2,y_1<y_2)+2I(x_2<x_1,y_2<y_1)-1,
}
the representation \eqref{urepres} yields  Kendall's $\tau$ at lag $k$ , that  is 
 
\al{
\tau_k ={}&  \Prob[(X_0^{(1)}-X_0^{(2)})(X_k^{(1)}-X_k^{(2)})>0]-\Prob[(X_0^{(1)}-X_0^{(2)})(X_k^{(1)}-X_k^{(2)})<0]\\
={}&2\Prob[X_0^{(1)}<X_0^{(2)},X_k^{(1)}<X_k^{(2)}]+2\Prob[X_0^{(2)}<X_0^{(1)},X_k^{(2)}<X_k^{(1)}]-1
}
The corresponding spectral density  will be denoted by 
\als{
\label{spectau}  
\spec_{\tau} (\omega ) = \frac{1}{2\pi}\sum_{k\in\Z} \tau_ke^{-ik\omega}.
}
As the distribution function $F$ and $F_k$ are assumed to be continuous, $\tau_k$ can also be represented in the form \eqref{urepres} using the kernel 
\als{\label{kerneltaualt} 
h\big ( ({x_1},{y_1})^T , ({x_2},{y_2})^T \big)=4\big(I(x_1<x_2)-\frac{1}{2}\big)\big(I(y_1<y_2)-\frac{1}{2}\big)
}
\item[(iii)] If $m=3$, $\Gamma\{i,j,k\}$ denotes the set of all permutations of $\{i,j,k\}$
and the kernel $h$ is defined by
\als{\label{hrho}
h\big ( ({x_1},{y_1})^T , ({x_2},{y_2})^T, ({x_3},{y_3})^T \big)
=\frac{1}{6}\sum_{\gamma\in\Gamma\{1,2,3\}}[12I(x_{\gamma(1)}<x_{\gamma(2)},y_{\gamma(1)}<y_{\gamma(3)})-3].
}
we obtain the (lag $k$) population version of Spearman's $\rho$, that 
is  
\[
\rho_k =3(\Prob[(X_0^{(1)}-X_0^{(2)})(X_k^{(1)}-X_k^{(3)})>0]-\Prob[(X_0^{(1)}-X_0^{(2)})(X_k^{(1)}-X_k^{(3)})<0]),
\]
The corresponding spectral density will be denoted by 
\als{
\label{specrho}  
\spec_{\rho} (\omega ) = \frac{1}{2\pi}\sum_{k\in\Z} \rho_ke^{-ik\omega}.
}
Given continuity of $F$ and $F_k$, $\rho_k$ can also be represented in the form \eqref{urepres} using the following kernel 
\als{\label{kernelrhoalt} 
h\big ( ({x_1},{y_1})^T , ({x_2},{y_2})^T, ({x_3},{y_3})^T \big)
=\sum_{\gamma\in\Gamma\{1,2,3\}}\big[2\big(I(x_{\gamma(1)}<x_{\gamma(2)})-\frac{1}{2}\big)\big(I(y_{\gamma(1)}<y_{\gamma(3)})-\frac{1}{2}\big)\big]
}
\ei
}
\end{example}

\medskip

In the remaining part of the manuscript we estimate the dependence measures $\xi_k$ (at lag $k$) by a $U$-statistic of order $m$, that is 
\als{ 
\xi_{n,k}={}&U_{n-|k|}(h)=\frac{1}{\binom{n-|k|}{m}}\sum\limits_{\substack{t_1,\dots,t_m\in\mathcal{T}_k\\t_1<\dots<t_m}}h\Big(\colvec{X_{t_1}}{X_{t_1+k}},\dots,\colvec{X_{t_m}}{X_{t_m+k}}\Big).
\label{ustatk}
}
Estimates of corresponding spectral densities are defined as in \eqref{per}. The asymptotic properties of such spectral density estimates are investigated in the following section. \\ 
Before proceeding, we recall the Hoeffding decomposition for U-statistics. Recall that $h$ is a symmetric kernel of order $m$ and let $\bm{Y}^{(1)},\dots,\,\bm{Y}^{(m)}$ denote independent identically distributed copies of $\colvec{X_0}{X_k}\sim F_k$. We now recursively define kernels $h_{c,k}$ by 
\begin{align} \label{hkhoef}
h_{c,k}(\bm{y_1},\dots,\bm{y_c}):={}&\E[h(\bm{Y}^{(1)},\dots,\bm{Y}^{(m)})|\bm{Y}^{(1)}=\bm{y_1},\dots,\bm{Y}^{(c)}=\bm{y_c}]\\
\nonumber 
&{}-\sum_{j=1}^{c-1}\sum\limits_{\substack{\{\nu_1,\dots,\nu_j\}\subset\{1,\dots,c\}\\ \nu_1<\dots<\nu_j}}h_{j,k}(\bm{y_{\nu_1}},\dots,\bm{y_{\nu_j}})-\xi_k\\ \nonumber 
={}&\int_{\R^2}\cdots\int_{\R^2}h(\bm{u_1},\dots,\bm{u_m})\prod_{j=1}^{c}(dG_{\bm{y_j}}(\bm{u_j})-dF_k(\bm{u_j}))\prod_{j=c+1}^{m}dF_k(\bm{u_j})~,
\end{align}
where $G_{\bm{y_j}}$ denotes the distribution of the Dirac measure  at $\bm{y_j}$.  If  
\begin{align*}
&U_{n-|k|}^{(c)}(h_{c,k})=\frac{1}{\binom{n-|k|}{c}}\sum_{\substack{t_1,\dots,t_c\in\mathcal{T}_k\\t_1<\dots <t_c}}h_{c,k}\Big(\colvec{X_{t_1}}{X_{t_1+k}},\dots,\colvec{X_{t_c}}{X_{t_c+k}}\Big)
\end{align*}
is the U-statistic based on the kernel $h_{c,k}$ we obtain for the statistic in \eqref{ustatk} the decomposition 
[see, for example  \cite{lee90}]
\als{\label{hdecomp}
\xi_{n,k}-\xi_k=\frac{m}{n-|k|}\sum_{t\in\mathcal{T}_k}h_{1,k}\colvec{X_t}{X_{t+k}}+\sum_{c=2}^m\binom{m}{c}U_{n-|k|}^{(c)}(h_{c,k}),
}
which will be an important tool in the  asymptotic analysis of the following sections.

\section{Asymptotic theory for U-lag-window estimates} \label{sec3}
\def\theequation{3.\arabic{equation}}
\setcounter{equation}{0}

\subsection{Consistency of U-lag-window estimates} \label{sec:cons}

Our first main result 
shows  that for a general class of symmetric kernels the statistic 
 $\hat{f}_{n,\xi}$ consistently estimates the spectral density $ \spec_{\xi} $ defined in \eqref{uspecdens} if the following assumptions are
satisfied. 

\bi
\item[(C1)] The lag window $w_n(\cdot)$ can be written in the form
$w_n(k)=w\big (\frac{k}{r_n}\big )$, where 
 $w(\cdot)$ is a  uniformly continuous function,  supported on the interval  $[-1,1]$,  satisfying
  $\|w\|_{\infty} \leq 1$, $w(0) = 1$, $w(-x)=w(x)$ for all $x\in\R$, and  
$r_n=n^{\frac{1}{2}-\nu}$ for some  $\nu \in (0,\frac{1}{2})$.
\item[(C2)] 
There exist  constants $\delta, M_0>0$ such that for all $t_1,\dots,t_{m},k\in\Z$, $1\leq j\leq 2m$,
\al{
\max\Big\{\int_{\R}\dots\int_{\R}|h|^{2+\delta}dG,
\int_{\R}\dots\int_{\R}|h|^{2+\delta}dG_j^{(1)}dG_j^{(2)}\Big\}\leq M_0<\infty,
}
where $G$, $G_j^{(1)}$ and $G_j^{(2)}$ denote the joint distributions of $(X_{t_{(1)}},\dots,X_{t_{(2m)}})$, $(X_{t_{(1)}},\dots,X_{t_{(j)}})$ and  $(X_{t_{(j+1)}},\dots,X_{t_{(2m)}})$, respectively,  and  $t_{(1)}\leq\dots\leq t_{(2m)}$ is  the order 
statistic  of $\{t_1,t_1+k,t_2,t_2+k\,\dots,t_m,t_m+k\}$.
\item[(C3)] The process $\proc$ is $\beta$-mixing and for some $\delta^{\prime}<\delta$ 
with $\beta$-mixing coefficients satisfying $\beta(n)=O(n^{-\gamma}),$ where $\gamma=\frac{2+\delta^{\prime}}{\delta^{\prime}}$.
\ei

\smallskip 

\begin{satz}\label{consistent}
If Assumptions (C0) -- (C3) are satisfied, we have  for any fixed $\omega\in\R$ 
\[
\uper\overset{\Prob}{\to} \spec_{\xi}(\omega),\qquad (n\rightarrow\infty).
\]
\end{satz}

\subsection{Asymptotic distribution of U-lag-window estimates} \label{sec:norm}

In this section we establish asymptotic normality of the spectral density estimators. Throughout this section we focus our attention on settings where $\xi$ is Kendall's $\tau$ or Spearman's $\rho$. Recalling  the discussion in Example \ref{exustat} it follows that $\tau_k $ und $\rho_k $
can be estimated by the $U$-statistics 
\al{
\tau_{n,k}\overset{a.s.}{=}\frac{1}{\binom{n-|k|}{2}}\underset{\substack{t_1,t_2\in\mathcal{T}_k\\ t_1<t_2}}{\sum}4\Big( I(X_{t_1}<X_{t_2})-\frac{1}{2}\Big)\Big(I(X_{t_1+k}<X_{t_2+k})-\frac{1}{2}\Big),
}
and
\al{
\rho_{n,k}\overset{a.s.}{=}\frac{1}{\binom{n-|k|}{2}}\underset{\substack{t_1,t_2,t_3\in\mathcal{T}_k\\ t_1<t_2<t_3}}{\sum}\sum_{\gamma\in\Gamma\{1,2,3\}}2\Big(I(X_{t_{\gamma(1)}}<X_{t_{\gamma(2)}})-\frac{1}{2}\Big)\Big(I(X_{t_{\gamma(1)}+k}<X_{t_{\gamma(3)}+k})-\frac{1}{2}\Big),
}
respectively. Note that these  $U$-statistics  have  bounded  kernels,  satisfy Assumption (C2) and can  be written as a product or sum of products of two centered functions of random variables. This special structure is crucial for obtaining the asymptotic distribution results given below. It is not clear if similar results hold without imposing this kind of structure on the kernel $h$.\\
Throughout this section we write $\xi_k$ if assumptions or results are the same for both Kendall's $\tau$ and Spearman's $\rho$. On the other hand we explicitly write $\tau$ or $\rho$ if the results or arguments are different. For example, from \eqref{hdecomp} we obtain for Kendall's $\tau$ and Spearman's $\rho$ the decomposition
\[
\hat{f}_{n,\xi}(\omega)=\frac{1}{2\pi}\sum_{|k|\leq \floor{r_n}} \wk\Big\{\xi_k+\frac{m}{n-|k|}\sum_{t\in\mathcal{T}_k}h_{1,k}^{\xi}\colvec{X_t}{X_{t+k}}+\sum_{c=2}^m\binom{m}{c}U_{n-|k|}^{(c)}(h_{c,k})\Big\}e^{-ik\omega},
\]
and therefore only $\xi_k$ appears in the formula. We will demonstrate that under suitable assumptions the term corresponding to the linear part converges to a normal distribution and that the term corresponding to the degenerate part is asymptotically negligible. 
In what follows we assume that (C0) -- (C3) hold and impose the following additional conditions. 

\bi
\item[(N1)] The process $\proc$ is $\alpha$-mixing  with corresponding $\alpha$-mixing coefficients satisfying 
$
\alpha(n)=O(n^{-\nu}),
$
where $\nu>7$. 
\item[(N2)] For the lag window generator $w$ there exists a 'characteristic exponent' $d>0$ being the largest integer such that 
\[
C_{w}(d):=\lim_{u\to 0}\frac{1-w(u)}{|u|^d}
\]
exists, is finite and non-zero.  For this $d$ we have $\sum_{k\in\Z}|k|^d|\xi_k|<\infty$. 
\item[(N3)] $r_n = o(n^\theta)$ where $\theta = \min\Big\{\frac{2(\delta-\delta^{\prime})}{\delta^{\prime}(2+\delta)},1\Big\}$ and $\delta, \delta'$ are from conditions (C2),(C3).
\ei

\medskip 

\begin{rem}\label{remonassN} {\rm 
The summability condition $\sum_{k\in\Z}|k|^d|\xi_k|<\infty$ in assumption (N2) implies the existence of the 'generalized $d^{th}$ derivative' of $\spec_{\xi}(\omega)$ 
\[
\spec_{\xi}^{[d]}(\omega):=\frac{1}{2\pi}\sum_{k\in\Z}|k|^d\xi_ke^{-ik\omega}
\]
and can thus be interpreted as a smoothness condition; note that for even $d$ this coincides with the usual $d$'th order derivative. The other part of assumption (N2) places mild restrictions on the lag-window generator for which the rate of the scale parameter is limited by assumption (N3). Note that (N3) is satisfied for scale parameters leading to optimal asymptotic mean squared error rates (see Remark \ref{MSE}).
}
\end{rem}

We begin by examining the asymptotic distribution of $\hat{f}_{n,\rho}(\omega)$.

\begin{satz}\label{asympnormalityrho}
Assume that conditions  (C0) -- (C3) and (N1) -- (N3) are satisfied and that $\omega\in(-\pi, \pi]$.
If  $\spec_{\rho}(\omega)\neq0$, then 
\be\label{weakconvrho}
\sqrt{\frac{n}{r_n}}\Big(\hat{f}_{n,\rho}(\omega)-\spec_{\rho}(\omega)-b_{\rho}(\omega)\Big)\overset{\mathcal{D}}{\longrightarrow}\mathcal{N}(0,\sigma_{\rho}^2(\omega)),
\ee
where 
\bea
\sigma_{\rho}^2(\omega)=(1+I(\omega\in\{0,
\pi\}))\spec_{\rho}^2(\omega)\int_{-1}^1w^2(x)dx.
\eea
and
\bea
b_{\rho}(\omega) := \E[\hat{f}_{n,\rho}(\omega)] - \spec_{\rho}(\omega) = -C_w(d)r_n^{-d}\spec^{[d]}_{\rho}(\omega)+o(r_n^{-d})=:r_n^{-d}b_{\rho,\omega}+o(r_n^{-d})
\eea
If $\spec_{\rho}(\omega)=0$, we have 
\al{
\sqrt{\frac{n}{r_n}}\Big(\hat{f}_{n,\rho}(\omega)-\E[\hat{f}_{n,\rho}(\omega)]\Big)\overset{\Prob}{\longrightarrow}0.
}
\end{satz}

Interestingly, the limiting distribution has exactly the same form as the limiting distribution for the usual spectral density where $\xi$ corresponds to covariance. This is remarkable, since Spearman's $\rho$ is based on covariances of \textit{ranks}. Asymptotic normality of $\hat{f}_{n,\rho}(\omega)$ was also obtained in~\cite{DetteEtAl2016} under a different set of assumptions on the serial dependence and using a completely different set of proof techniques. Specifically, their results require dependence to decay exponentially. The next result establishes asymptotic normality of $\hat{f}_{n,\xi}(\omega)$ with $\xi$ corresponding to Kendall's $\tau$. The asymptotic distribution of $\hat{f}_{n,\tau}(\omega)$ cannot be obtained from the findings in~\cite{DetteEtAl2016} (under any assumptions) since Kendall's $\tau$ is a non-linear functional of the copula.   
  
\begin{satz}\label{asympnormalitytau}
Assume that conditions (C0) -- (C3) and (N1) -- (N3) are satisfied and that $\omega\in(-\pi, \pi]$. If $\spec_{\rho}(\omega)\neq0$, then 
\be\label{weakconv}
\sqrt{\frac{n}{r_n}}\Big(\hat{f}_{n,\tau}(\omega)-\spec_{\tau}(\omega)-b_{\tau}(\omega)\Big)\overset{\mathcal{D}}{\longrightarrow}\mathcal{N}(0,\sigma_{\tau}^2(\omega)),
\ee
where
\bea
\sigma_{\tau}^2(\omega) = \frac{4}{9}(1+I(\omega\in\{0,
\pi\}))\spec_{\rho}^2(\omega)\int_{-1}^1w^2(x)dx
\eea
and
\bea
b_{\tau}(\omega) := \E[\hat{f}_{n,\tau}(\omega)] - \spec_{\tau}(\omega) = -C_w(d)r_n^{-d}\spec^{[d]}_{\tau}(\omega)+o(r_n^{-d})=:r_n^{-d}b_{\tau,\omega}+o(r_n^{-d}).
\eea
If $\spec_{\rho}(\omega)=0$, we have 
\al{
\sqrt{\frac{n}{r_n}}\Big(\hat{f}_{n,\tau}(\omega)-\E[\hat{f}_{n,\tau}(\omega)]\Big)\overset{\Prob}{\longrightarrow}0.
}
\end{satz}

It is remarkable that the asymptotic variance of estimator $\hat{f}_{n,\tau}(\omega)$ depends on the spectral measure  $\spec_{\rho}(\omega) $ obtained from Spearman's $\rho$ (provided that $\spec_{\rho}(\omega) \neq 0$). This is in sharp contrast to the finding in Theorem~\ref{asympnormalityrho} and spectral density estimation based on covariances. The results in Theorem~\ref{asympnormalitytau} provide an asymptotic analysis of the estimator introduced in~\cite{ahdesmaki2005}. We conclude this section by commenting on the optimal choice of window length $r_n$.

\begin{rem}\label{MSE} {\rm  ~ Both theorems allow to determine the scale parameter $r_n$ such that the asymptotic mean squared error is minimized. To be precise, define $\sigma_{\xi,\omega}^2:=\sigma_{\xi}^2(\omega)$. Then the asymptotic mean squared error takes the form  
\al{
[r_n^{-2d}b_{\xi,\omega}^2+\frac{r_n}{n}\sigma_{\xi,\omega}^2](1+o(1)).
} 
Assuming that $b_{\xi,\omega}\neq 0$, we obtain that this expression is minimized for
\al{
r_n=\Big(\frac{2db_{\xi,\omega}^2}{\sigma_{\xi,\omega}^2}n\Big)^{\frac{1}{2d+1}}.
}
Note that for $d=2$ the asymptotic MSE is of the order $n^{-4/5}$. In that case the above scale parameter $r_n$ is of order $n^{1/5}$ and satisfies Assumptions (C1) and (N3) if the mixing coefficients $\beta(n)$ decay sufficiently quickly. More precisely, as for Kendall's $\tau$ and Spearman's $\rho$ the kernels of the $U$-statistic are bounded, we can choose $\delta$ in Assumption (C2) arbitrarily large. Assumption (N3) is satisfied if $1/5 < \theta = \min\Big\{\frac{2(\delta-\delta^{\prime})}{\delta^{\prime}(2+\delta)},1\Big\}$, which is equivalent to $\delta^{\prime} < 10 \delta/(12 + \delta)$. Since $\delta$ is arbitrarily large, we can choose any $\delta^{\prime} < 10$ and (N3), (C3) will hold if $\beta(n) = O(n^{-\gamma})$ for some $\gamma > 6/5$. 
}
\end{rem}

\section{Proofs} \label{sec4}
\def\theequation{4.\arabic{equation}}
\setcounter{equation}{0}

\subsection{Proof of Theorem \ref{consistent}}

We first  illustrate the main steps in the proofs. These rely on several delicate bounds,
which will be shown below. Rearranging sums in (\ref{per}) and using assumption (C1), the U-lag-window estimate can be decomposed as follows
\bea
\uper -\spec_{\xi}(\omega) 
&=& s_{n,1} - s_{n,2} +s_{n,3}
\eea
where  
\bea
s_{n,1} & = &   \frac{1}{2\pi}\sum_{|k|<n}\Big (w\Big (\frac{k}{r_n}\Big)-1\Big )\xi_ke^{-ik\omega}, \\
s_{n,2} & = & \frac{1}{2\pi}\sum_{|k|\geq n}\xi_ke^{-ik\omega}, \\
s_{n,3} & = & \frac{1}{2\pi}\sum_{|k|\leq \floor{r_n}}w\Big  (\frac{k}{r_n}\Big )(\xi_{n,k}-\xi_k)e^{-ik\omega}.
\eea
We will show in Section \ref{sn12sec} that 
\beq\label{sn12}
s_{n,1}\overset{a.s.}{\to} 0\qquad\text{and}\qquad s_{n,2}\overset{a.s.}{\to} 0. 
\eeq
For a proof of   $
s_{n,3}\overset{\Prob}{\to}0$ 
we use the  Hoeffding decomposition \eqref{hdecomp},
which gives 
\als{
2\pi s_{n,3}= & \sum_{|k|<n}w\Big(\frac{k}{r_n}\Big)(\xi_{n,k}-\xi_k)e^{-ik\omega}= d_{n,1}+d_{n,2}.\label{convd}
}
where
\als{
d_{n,1} ={}&  \sum_{|k|\leq \floor{r_n}} \wk \frac{m}{n-|k|}\sum_{t\in\mathcal{T}_k}h_{1,k}\colvec{X_t}{X_{t+k}}e^{-ik\omega}\notag\\
d_{n,2} ={}& \sum_{|k|\leq \floor{r_n}} \wk \sum_{c=2}^m\binom{m}{c}U_{n-|k|}^{(c)}(h_{c,k})e^{-ik\omega}\notag
}
The assertion of the theorem now follows from the estimates
\als{\label{dn2}
d_{n,2} &=O_{\Prob}(r_nn^{-1/2-\theta/2})=o_{\Prob}(1),\\
\label{dnlin}
d_{n,1} &=O_{\Prob}(r_nn^{-1/2})=o_{\Prob}(1),
}
which are shown in Section \ref{dn2sec} and  \ref{dnlinsec}, respectively.





\subsubsection{Proof of (\ref{sn12})} \label{sn12sec} 
Using the fact that $w(0)=1$ and $\sup_{|k|<n}\big(\big|w\big(\frac{k}{r_n}\big)\big|+1\big)\leq 2\quad  \alle n\in\N$, we obtain for any fixed $0\leq K<n$
\al{
|s_{n,1}| \leq{}& 
\sum_{|k|<n}\Big|w\Big(\frac{k}{r_n}\Big)-1\Big||\xi_k|
\leq{} \sum_{|k|\leq K}\Big|w\Big(\frac{k}{r_n}\Big)-w(0)\Big|\sup_{|k|\leq K}|\xi_k|+2\sum_{n>|k|>K}|\xi_k|\notag\\
\leq{} &(2K+1)\sup_{|k|\leq K}\Big|w\Big(\frac{k}{r_n}\Big)-w(0)\Big|\sup_{|k|\leq K}|\xi_k|+2\sum_{n>|k|>K}|\xi_k|\notag
}
As the lag window generator $w(\cdot)$ is continuous at $0$ we obtain for any fixed $K$
\als{\label{ersterlimes}
\lim_{n\rightarrow\infty}\Big [(2K+1)\sup_{|k|\leq K}\Big |w\Big (\frac{k}{r_n}\Big )-w(0)\Big 
|\sup_{|k|\leq K}|\xi_k|+2\sum_{n>|k|>K}|\xi_k|\Big]=0+2\sum_{|k|>K}|\xi_k|.
}
As inequality (\ref{ersterlimes}) holds for all $K$, we can conclude that $s_{n,1}\overset{a.s.}{\to}0$, as the $\xi_k$'s are absolutely summable. By the same argument, we also have
$
|s_{n,2}|  \leq \frac{1}{2\pi}\sum_{|k|\geq n}|\xi_k|\to 0
$
for $n\rightarrow\infty$.

\subsubsection{Proof of (\ref{dn2})}\label{dn2sec}

The proof  is based on an extension to lagged data of a covariance inequality by \cite{yoshihara76}. 
More precisely, we prove in the technical Appendix, Section~\ref{sec:proofyoshi}, that for fixed $2\leq c\leq m$
\als{\label{laggedyoshihara}
\sup_{|k|\leq \floor{r_n}}\E\Big[\Big(U_{n-|k|}^{(c)}(h_{c,k})\Big)^2\Big]=O(n^{-1-\theta}).
} 
where $\theta=\min\Big\{\frac{2(\delta-\delta^{\prime})}{\delta^{\prime}(2+\delta)},1\Big\}$. Note that the above bound holds uniformly over a growing number of lags $k$ while the result in~\cite{yoshihara76} only holds for a \textit{fixed} $k$. Observe that,
\al{
\E|d_{n,2}|
\leq& (2r_n+1)\sum_{c=2}^m\binom{m}{c} \sup_{|k|\leq\floor{r_n}}\E|U_{n-|k|}^{(c)}(h_{c,k})| 
\leq (2r_n+1) C_m
 \sup_{|k| \leq \floor{r_n}} \big (\E [(U_{n-|k|}^{(c)}(h_{c,k}) )^2 ]\big )^{\frac{1}{2}}
}
where the constant $C_m$ does not depend on $m$.
Consequently, equation \eqref{laggedyoshihara} yields
$
\E|d_{n,2}|=O(r_nn^{-1/2-\theta/2}),
$
which establishes (\ref{dn2}).

\medskip

\subsubsection{Proof of (\ref{dnlin})} \label{dnlinsec} 
Introduce the notation $\bm{X}_{t,k}:=(X_{t},X_{t+k})^T$. We only consider positive lags $k$, negative $k$ can be treated analogously. Similar arguments as in the proof of (\ref{dn2}) yield
\al{
&\E\Big|\sum_{\kip} \wk \frac{m}{n-|k|}\sum_{t\in\mathcal{T}_k}h_{1,k}(\bm{X}_{t,k})e^{-ik\omega}\Big|
\\
\leq{}& 
(r_n+1)\sup_{\kip}\Big(\E\Big[\Big(\frac{m}{n-k}\sum_{t=1}^{n-k}h_{1,k}(\bm{X}_{t,k})\Big)^2\Big]\Big)^{\frac{1}{2}}.
}
Next,
\al{
&\E\Big[\Big(\frac{m}{n-k}\sum_{t=1}^{n-k}h_{1,k}(\bm{X}_{t,k})\Big)^2\Big]\notag\\
={}&\frac{m^2}{(n-k)^2}\sum_{t=1}^{n-k}\E\Big[\Big(h_{1,k}(\bm{X}_{t,k})\Big)^2\Big]+\frac{2m^2}{(n-k)^2}\sum_{l=1}^{n-k-1}\sum_{u=1}^{n-k-l}\E\Big[h_{1,k}(\bm{X}_{u,k})h_{1,k}(\bm{X}_{u+l,k})\Big]
}
As $(\E|Z|^p)^{\frac{1}{p}}\leq (\E|Z|^q)^{\frac{1}{q}}$ for $p<q$, we have
\begin{eqnarray*}
\big|\E\big[h_{1,k}(\bm{X}_{t,k})\big]\big|^2&\leq& \big(\E\big|h_{1,k}(\bm{X}_{t,k})\big|^{2}\big)^{\frac{1}{2}}\big(\E\big|h_{1,k}(\bm{X}_{t,k})\big|^{2}\big)^{\frac{1}{2}} \\
&\leq & \big(\E\big|h_{1,k}(\bm{X}_{t,k})\big|^{2+\delta}\big)^{\frac{1}{2+\delta}}\big(\E\Big|h_{1,k}(\bm{X}_{t,k})\big|^{2+\delta}\big)^{\frac{1}{2+\delta}}
\leq M_1^{\frac{2}{2+\delta}},
\end{eqnarray*}
which gives
\begin{equation}
\sup_{\kip}\frac{m^2}{(n-k)^2}\sum_{t=1}^{n-k}\E\big[h_{1,k}(\bm{X}_{t,k})\big]^2\leq \sup_{\kip}\frac{m^2}{n-k}M_1^{\frac{2}{2+\delta}}=O(n^{-1}).\label{var}
\end{equation}
The following bound will be established in Lemma~\ref{boundshm} in the Appendix (see Section~\ref{sec:proofboundshm})
\[
\E\Big|h_{1,k}(\bm{X}_{u,k})h_{1,k}(\bm{X}_{u+l,k})\Big|\leq \begin{cases}
2M_1^{\frac{2}{2+\delta}}\beta^{\frac{\delta}{2+\delta}}(l-k),& \text{if } l>k\geq 0,\\
8M_1^{\frac{2}{2+\delta}}\beta^{\frac{\delta}{2+\delta}}(\min\{l,k-l\}),& \text{if } 0\leq l\leq k.
\end{cases}
\]
Thus
\al{
&\Big|\frac{2m^2}{(n-k)^2}\sum_{l=1}^{n-k-1}\sum_{u=1}^{n-k-l}\E\big[h_{1,k}(\bm{X}_{u,k})h_{1,k}(\bm{X}_{u+l,k})\big]\Big|
\\
\leq{}&\hspace{0cm} \frac{2m^2}{(n-k)^2}\sum_{l=1}^{k}\sum_{u=1}^{n-k-l}\E\big|h_{1,k}(\bm{X}_{u,k})h_{1,k}(\bm{X}_{u+l,k})\big|
\\
&\hspace{1cm}+\frac{2m^2}{(n-k)^2}\sum_{l=k+1}^{n-k-1}\sum_{u=1}^{n-k-l}\E\big|h_{1,k}(\bm{X}_{u,k})h_{1,k}(\bm{X}_{u+l,k})\big|
\\
\leq{}&\hspace{0cm} \frac{16m^2}{(n-k)^2}M_1^{\frac{2}{2+\delta}}\sum_{l=1}^{k}\sum_{u=1}^{n-k-l}\beta^{\frac{\delta}{2+\delta}}(\min\{l,k-l\})
+\frac{4m^2}{(n-k)^2}M_1^{\frac{2}{2+\delta}}\sum_{l=k+1}^{n-k-1}\sum_{u=1}^{n-k-l}\beta^{\frac{\delta}{2+\delta}}(l-k)
\\
\leq{}&\hspace{0cm} \frac{16m^2}{(n-k)^2}M_1^{\frac{2}{2+\delta}}(n-k)2\sum_{v=0}^{\floor{\frac{k}{2}}}\beta^{\frac{\delta}{2+\delta}}(v)+\frac{4m^2}{(n-k)^2}M_1^{\frac{2}{2+\delta}}(n-k)\sum_{v=1}^{n-2k-1}\beta^{\frac{\delta}{2+\delta}}(v).
}
By assumption (C3) $\sum_{j=1}^{\infty}\beta^{\frac{\delta}{2+\delta}}(j)<\infty$. Therefore, we have
\begin{equation}
\sup_{\kip}\frac{2m^2}{(n-k)^2}\sum_{l=1}^{n-k-1}\sum_{u=1}^{n-k-l}\E\big[h_{1,k}(\bm{X}_{u,k})h_{1,k}(\bm{X}_{u+l,k})\big]=O(n^{-1}).\label{cov}
\end{equation}
Equations (\ref{var}) and (\ref{cov}) yield
$
\E|d_{n,1}|=O(r_nn^{-\frac{1}{2}})
$
and the assertion follows observing that  $r_n=n^{\frac{1}{2}-\nu}$. \hfill $\Box$


\subsection{Proof of Theorem \ref{asympnormalitytau} and \ref{asympnormalityrho} - main arguments}

In the following proof we write $\xi$ if the results hold for general dependence measures that fulfill the assumptions (C0) -- (C3) and (N1) -- (N3).  Otherwise we explicitly write $\tau$ or $\rho$.\\
Under Assumption (N3) and with (\ref{dn2}),
\al{
\sqrt{\frac{n}{r_n}}\Big(\frac{1}{2\pi}\sum_{|k|\leq \floor{r_n}}\wk\sum_{c=2}^m\binom{m}{c}U_{n-|k|}^{(c)}(h_{c,k}^{\xi})e^{-ik\omega}\Big)\overset{\Prob}{\longrightarrow}0 \quad (n\rightarrow\infty).
}
Furthermore, in Section~\ref{sec:prooftriangarray} we will prove that
\als{\label{triangarray}
\hat{f}_{n,\xi}(\omega)=\tilde{f}_{n,\xi}(\omega)+o_{\Prob}\Big(\sqrt{\frac{r_n}{n}}\Big),
}

where $\tilde{f}_{n,\xi}(\omega)=\frac{1}{2\pi}\sum_{|k|\leq r_n}\wk \Big\{\xi_k+\frac{m}{n}\sum_{t=1}^nh_{1,k}^{\xi}(\bm{X}_{t,k})\Big\}e^{-ik\omega}$. Then, 

\al{
\tilde{f}_{n,\xi}(\omega)-\E[\tilde{f}_{n,\xi}(\omega)]={}&\frac{1}{2\pi}\sum_{|k|\leq r_n}\wk\frac{m}{n}\sum_{t=1}^nh_{1,k}^{\xi}(\bm{X}_{t,k})e^{-ik\omega}
=:{}\sum_{t=1}^nW_{n,t}^{\xi}(\omega),
}
where, by construction, the random variables $(W_{n,t}^{\xi})_{t=1,\dots,n}$ form a triangular array of $\beta$-mixing random variables with mixing coefficients $\beta^W(u)\leq\beta^X(0\vee (u-2r_n))$. To prove the asymptotic normality, we will apply the blocking technique described in Section \ref{blocking}. That is, we choose $\mu_n$ blocks of length $p_n$ and $\mu_n$ blocks of length $q_n$ such that
\bea
q_n/p_n\rightarrow 0,\qquad r_n/q_n\rightarrow 0,\qquad p_n/n\rightarrow 0,\qquad \mu_n\beta(q_n)\rightarrow 0.
\eea
According to Assumptions (C1), (C3) and (N3) one possible choice is $r_n=O(n^{1/2-\nu}),\,0<\nu<\min\{\theta,\frac{1}{2}\}$, $q_n=O(n^{1/2})$, $p_n=O(n^{1/2+\nu})$. Then we decompose
\als{\label{blocksum}
\sum_{t=1}^nW_{n,t}^{\xi}(\omega)=\sum_{j=1}^{\mu_n}\sum_{t\in\Gamma_j}W_{n,t}^{\xi}(\omega)+\sum_{j=1}^{\mu_n}\sum_{t\in\Delta_j}W_{n,t}^{\xi}(\omega)+\sum_{t\in\mathcal{R}}W_{n,t}^{\xi}(\omega).
}
Next, we show that the remaining part and the part corresponding to the \az{small} blocks are negligible whereas the \az{big} blocks satisfy the Lyapunov condition and yield the asymptotic variance. Observe that
\als{\label{symmetryh1}
h_{1,k}^{\tau}(\bm{X}_{t,k})\overset{\mathcal{D}}{=}h_{1,-k}^{\tau}(\bm{X}_{t,-k})\text{ and }h_{1,-k}^{\rho}(\bm{X}_{t,k})\overset{\mathcal{D}}{=}h_{1,-k}^{\rho}(\bm{X}_{t,-k})
}
for all $k,t\in\Z$, and hence, $\sum_{t=1}^nW_{n,t}^{\xi}(\omega)$ is real and symmetric in $\omega$. To prove~\eqref{symmetryh1} observe that by stationarity of $\{X_t\}_{t\in \Z}$ we have that $(\bm{X}_{0,k})\eid(\bm{X}_{0,-k})$ and $(\bm{X}_{0,k}^{(1)})\eid(\bm{X}_{0,-k}^{(1)})$ and hence, $ \tau_k=\tau_{-k}$, $
\rho_k=\rho_{-k}$.
Consequently, we obtain in the case of Kendall's $\tau$ and Spearmans's $\rho$ 
$h_{1,k} (({x},{y})^T) =h_{1,-k} (({y},{x})^T) $, which yields  
\bea
h_{1,k}  (\bm{X}_{t,k}) \eid h_{1,k} (\bm{X}_{t-k,k}) 
\eid h_{1,-k}  (\bm{X}_{t,-k})
\eea

Moreover, we will show in section \ref{asymnormydet} that for any $a_n\rightarrow\infty$ with $a_n/n=o(1)$, $r_n/a_n=o(1)$ and $\omega\in(-\pi, \pi]$,
\als{\label{varblocktau}
\Big|\E\Big[\sum_{t_1=1}^{a_n}W_{n,t_1}^{\tau}(\omega)&\sum_{t_2=1}^{a_n}W_{n,t_2}^{\tau}(\omega)\Big]-\frac{a_nr_n}{n^2}\sigma_{\tau}^2(\omega)\Big|=o\Big(\frac{a_nr_n}{n^2}\Big),
}
\als{\label{varblockrho}
\Big|\E\Big[\sum_{t_1=1}^{a_n}W_{n,t_1}^{\rho}(\omega)&\sum_{t_2=1}^{a_n}W_{n,t_2}^{\rho}(\omega)\Big]-\frac{a_nr_n}{n^2}\sigma_{\rho}^2(\omega)\Big|=o\Big(\frac{a_nr_n}{n^2}\Big).
}
Next, the last summand in (\ref{blocksum}) contains at most $O(p_n+q_n)$ summands. Hence, by (\ref{varblocktau}) and (\ref{varblockrho}) we have,
\bea
\Var\Big(\sum_{t\in\mathcal{R}}W_{n,t}^{\xi}(\omega)\Big)=O\Big(\frac{(p_n+q_n)r_n}{n^2}\Big)=o\Big(\frac{r_n}{n}\Big),
\eea
that is $\sum_{t\in\mathcal{R}}W_{n,t}^{\xi}(\omega)=o_{\Prob}\Big(\sqrt{\frac{r_n}{n}}\Big)$. Next, we show that the sum over the small blocks is negligible. By Lemma \ref{link} with the function $\tilde g(\cdot)=I(\cdot\geq \varepsilon)$ we obtain
\bea
\Big|\Prob\Big(\sqrt{\frac{n}{r_n}}\sum_{j=1}^{\mu_n}\sum_{t\in\Delta_j}W_{n,t}^{\xi}(\omega)\geq\varepsilon\Big)- \Prob\Big(\sqrt{\frac{n}{r_n}}\sum_{j=1}^{\mu_n}\sum_{t\in\Delta_j}\zeta^{\xi}_{n,t}(\omega)\geq\varepsilon\Big)\Big|\leq(\mu_n-1)\beta^W(p_n),
\eea
where $\zeta^{\xi}_{n,t}(\omega)$ denote the random variables of the independent block sequence corresponding to the $\Delta$-blocks. By the assumptions on $p_n$ and $\beta^X$ the term on the right hand side in the above expression converges to 0. Observing that the variables $\zeta^{\xi}_{n,t}(\omega)$ are centered and (\ref{varblocktau}) or respectively (\ref{varblockrho}) applied to the independent blocks $\sum_{t\in\Delta_j}\zeta^{\xi}_{n,t}(\omega)$ yields
\bea
\Var\Big(\sqrt{\frac{n}{r_n}}\sum_{j=1}^{\mu_n}\sum_{t\in\Delta_j}\zeta^{\xi}_{n,t}(\omega)\Big)=\frac{n}{r_n}\sum_{j=1}^{\mu_n}\Var\Big(\sum_{t\in\Delta_j}W_{n,t}^{\xi}(\omega)\Big)=O\Big(\frac{\mu_nq_n}{n}\Big)=o(1),
\eea
where we have used the definition of $\zeta^{\xi}_{n,t}$ and the assumption that $q_n/p_n=o(1)$. Hence, it remains to prove that $\sqrt{\frac{n}{r_n}}\sum_{j=1}^{\mu_n}\sum_{t\in\Gamma_j}W_{n,t}^{\xi}(\omega)$ converges weakly. Note that for any measurable set $A$, by Lemma \ref{link} with function $g(\cdot)=I(\cdot\in A)$ and the assumptions on $q_n$ and $\beta^X$, we have
\bea
\Big|\Prob\Big(\sqrt{\frac{n}{r_n}}\sum_{j=1}^{\mu_n}\sum_{t\in\Gamma_j}W_{n,t}^{\xi}(\omega)\in A\Big)- \Prob\Big(\sqrt{\frac{n}{r_n}}\sum_{j=1}^{\mu_n}\sum_{t\in\Gamma_j}\zeta^{\xi}_{n,t}(\omega)\in A\Big)\Big|=o(1)
\eea
In order to prove the convergence in distribution of $\sqrt{\frac{n}{r_n}}\sum_{j=1}^{\mu_n}\sum_{t\in\Gamma_j}W_{n,t}^{\xi}(\omega)$, it suffices to show that the triangular array of independent random variables 
\bea
\Big(\sqrt{\frac{n}{r_n}}\sum_{t\in\Gamma_j}\zeta^{\xi}_{n,t}(\omega)\Big)_{j=1,\dots\mu_n}
\eea
satisfies the Lyapunov condition. To achieve this, we show that together with (\ref{varblocktau}) or (\ref{varblockrho}), respectively,
\als{\label{fourthmoment}
\sum_{j=1}^{\mu_n}\E\Big[\Big(\sum_{t\in\Gamma_j}\zeta^{\xi}_{n,t}(\omega)\Big)^4\Big]=O\Big(\frac{\mu_np_n^2r_n^2}{n^4}\Big).
}
and if $\spec_{\rho}(\omega)\neq 0$,
\als{\label{lowerboundvar}
\Big(\sum_{j=1}^{\mu_n}\Var\Big(\sum_{t\in\Gamma_j}\zeta^{\xi}_{n,t}(\omega)\Big)\Big)^2\geq c\frac{\mu_n^2p_n^2r_n^2}{n^4}
}
for some constant $c>0$ and $n$ sufficiently large. Hence, the Lyapunov condition is satisfied as $\mu_n\rightarrow \infty$ and we can conclude that the distribution of $\sqrt{\frac{n}{r_n}}\sum_{j=1}^{\mu_n}\sum_{t\in\Gamma_j}\zeta^{\xi}_{n,t}(\omega)$, where $\xi$ is either Kendall's $\tau$ or Spearman's $\rho$, converges weakly to a normal distribution, i.e.
\al{
\sqrt{\frac{n}{r_n}}\Big(\tilde{f}_{n,\xi}-\E[\tilde{f}_{n,\xi}]\Big)\overset{\mathcal{D}}{\to}\mathcal{N}(0,\sigma_{\xi}^2(\omega)).
}
If $\spec_{\rho}(\omega)=0$, it follows from equations (\ref{varblocktau}) and (\ref{varblockrho}) that 
\al{
\Var\Big(\sqrt{\frac{n}{r_n}}\Big(\tilde{f}_{n,\xi}-\E[\tilde{f}_{n,\xi}]\Big)\Big)=o(1)
}
and hence, by (\ref{triangarray}),
\al{
\sqrt{\frac{n}{r_n}}\Big(\hat{f}_{n,\xi}-\E[\hat{f}_{n,\xi}]\Big)\overset{\Prob}{\to} 0.
}
Finally, we have, for $\xi$ representing either $\tau$ or $\rho$,
\al{
\E[\tilde{f}_{n,\xi}(\omega)]={}&\frac{1}{2\pi}\sum_{|k|\leq r_n}\Big(\wk-1\Big)\xi_ke^{-ik\omega}+\frac{1}{2\pi}\sum_{|k|\leq r_n}\xi_ke^{-ik\omega}\\
={}&\spec_{\xi}(\omega)-\frac{1}{2\pi}\sum_{|k|> r_n}\xi_ke^{-ik\omega}+\frac{1}{2\pi}\sum_{|k|\leq r_n}\Big(\wk-1\Big)\xi_ke^{-ik\omega}.
}
Hence, the bias is given by
\al{
b_{\xi}(\omega)=\E[\tilde{f}_{n,\xi}(\omega)]-\spec_{\xi}(\omega)=\frac{1}{2\pi}\sum_{|k|\leq r_n}\Big(\wk-1\Big)\xi_ke^{-ik\omega}-\frac{1}{2\pi}\sum_{|k|> r_n}\xi_ke^{-ik\omega}.
}
Next, choose some $L_n\rightarrow\infty$ such that $L_n/r_n\rightarrow 0$. Then,
\als{\label{biais}
&\frac{1}{2\pi}\sum_{|k|\leq r_n}\Big(\wk-1\Big)\xi_ke^{-ik\omega}-\frac{1}{2\pi}\sum_{|k|> r_n}\wk\xi_ke^{-ik\omega}\notag\\
={}&r_n^{-d}\frac{1}{2\pi}\sum_{|k|\leq L_n}\frac{w(k/r_n)-1}{|k/r_n|^{d}}|k|^d\xi_ke^{-ik\omega}+r_n^{-d}\frac{1}{2\pi}\sum_{r_n\geq |k|> L_n}\frac{w(k/r_n)-1}{|k/r_n|^{d}}|k|^d\xi_ke^{-ik\omega}\notag\\
&-r_n^{-d}\frac{1}{2\pi}\sum_{|k|> r_n}\frac{r_n^d}{|k|^d}|k|^d\xi_ke^{-ik\omega}=:(I)+(II)+(III).
}
By Assumption (N2), $\sum_{k\in\Z}|k|^{q}|\xi_k|<\infty$ for $q\leq d$ and therefore,
\al{
|(III)|\leq O(r_n^{-d})\sum_{|k|>r_n}|k|^d|\xi_k|=o(r_n^{-d}).
}
For the second term in (\ref{biais}) as by Assumption (N2), $\sum_{k\in\Z}|k|^{q}|\xi_k|$ is finite for $q\leq d$ we obtain
\al{
|(II)|\leq O(1)r_n^{-d}\sup_{v\in [0,1]}\frac{w(v)-1}{|v|^d}\sum_{r_n\geq |k|> L_n}|k|^d|\xi_k|=o(r_n^{-d}),
}
where $\sup_{v\in [0,1]}\frac{w(v)-1}{|v|^d}$ is bounded since $w$ is bounded and the limit for $|v|\rightarrow 0$ exists by Assumption (N2). Finally,
\al{
(I)+r_n^{-d}C_w(d)\spec^{[d]}_{\xi}(\omega)={}&r_n^{-d}\frac{1}{2\pi}\sum_{|k|\leq L_n}\Big(\frac{w(k/r_n)-1}{|k/r_n|^{d}}+C_w(d)\Big)|k|^d\xi_ke^{-ik\omega}\\
&\quad+r_n^{-d}C_w(d)\frac{1}{2\pi}\sum_{|k|>L_n}|k|^de^{-ik\omega}\xi_k,
}
where the first summand is of order $o(r_n^{-d})$ since $\frac{L_n}{r_n}\rightarrow 0$ and $|k|^d\xi_k$ is absolutely summable. The second summand is of order $o(r_n^{-d})$ as $\sum_{k\in\Z}|k|^{q}|\xi_k|$ is finite. Hence,
\al{
b_{\xi}(\omega)=-r_n^{-d}C_w(d)\spec^{[d]}_{\xi}(\omega)+o(r_n^{-d}).
}
Conclude applying Slutsky's theorem. It remains to prove (\ref{varblocktau}), (\ref{varblockrho}), (\ref{fourthmoment}) and (\ref{lowerboundvar}). Detailed proofs of these results are given in the remaining part of this section.  \hfill $\Box$

\newpage

\subsection{Proof of (\ref{varblocktau}) -- (\ref{lowerboundvar})} \label{asymnormydet} 

The proofs of Theorem (\ref{varblocktau}) -- (\ref{lowerboundvar}) rely on two auxiliary results. These will be stated in this section whereas their detailed proof is deferred to the Appendix. The first Lemma bounds cumulants through $\alpha$-mixing coefficients, see Section~\ref{sec:proofcum} for a proof.

\begin{lem}\label{cum}
For $q\in\N$, let $(X_t^{(1)})_{t\in\Z},\dots,(X_t^{(q)})_{t\in\Z}$ be independent copies of a strictly stationary polynomially $\alpha$-mixing process $(X_t)_{t\in\Z}$ that are independent. For any $t\in\Z$, let $V_t:=(X_t,X_t^{(1)},\dots,X_t^{(q)})$. Then, for $p\in\N$, $t_1,\dots,t_p$ and measurable sets $A_1,\dots,A_p\subset\R^{q+1}$, there exists a constant $C_{p,q}$ such that 
\[
|\cum(I(V_{t_1}\in A_1),\dots,I(V_{t_p}\in A_p)|\leq C_{p,q}\alpha^{X}(\max_{i,j=1,\dots,p}|t_j-t_i|).
\]
\end{lem}

The Lemma that follows is a key observation which makes it possible to use theory from classical spectral density estimation in the case where the kernel $h$ can be written as a sum of a product of centered functions of random variables. This is a crucial insight for proving asymptotic normality of the estimators $\hat f_{n,\xi}$.

\begin{lem}\label{auxlemmavar}
Let $h$ denote a $U$-statistic of order $m$ and assume that $(X_t^{(1,j)})_{t\in\Z},\dots,(X_t^{(m-1,j)})_{t\in\Z}$, $j=1,\dots,q$ are independent copies of a strictly stationary process $(X_t)_{t\in\Z}$ that are mutually independent. Then, for $t_j,k_j\in\Z$, $j=1,\dots,q$, 

\als{\label{eprodhm}
\E\biggl[\prod_{j=1}^q\hunind{k_j}{X_{t_j}}{X_{t_j+k_j}}\biggr]={}&\E\bigg[\prod_{j=1}^{q}\biggl(h\Big(\binom{X_{t_j}}{X_{t_j+k_j}},\binom{X_{t_j}^{(1,j)}}{X_{t_j+k_j}^{(1,j)}},\dots,\binom{X_{t_j}^{(m-1,j)}}{X_{t_j+k_j}^{(m-1,j)}}\Big)-\xi_{k_j}\biggr)\biggr],\notag\\
&{}
}

where 
\[
\hunind{k_j}{X_{t_j}}{X_{t_j+k_j}}=\E\bigg[h\Big(\binom{X_{t_j}}{X_{t_j+k_j}},\binom{X_{t_j}^{(1,j)}}{X_{t_j+k_j}^{(1,j)}},\dots,\binom{X_{t_j}^{(m-1,j)}}{X_{t_j+k_j}^{(m-1,j)}}\Big)\biggl|\binom{X_{t_j}}{X_{t_j+k_j}}\biggr]-\xi_{k_j},
\] 
\[
\xi_{k_j}=\E\bigg[h\Big(\binom{X_{t_j}}{X_{t_j+k_j}},\binom{X_{t_j}^{(1,j)}}{X_{t_j+k_j}^{(1,j)}},\dots,\binom{X_{t_j}^{(m-1,j)}}{X_{t_j+k_j}^{(m-1,j)}}\Big)\biggr].
\]

In particular, if $(X_t^{(j)})_{t\in\Z}$, $j=1,\dots,5$ are independent copies of the strictly stationary process $\proc$ that are independent of each other,

\bi
\item[(i)] for Kendall's $\tau$
\al{
&\Cov\biggl(\hunindt{k_1}{X_{t_1}}{X_{t_1+k_1}},\hunindt{k_2}{X_{t_2}}{X_{t_2+k_2}}\biggr)\\
={}&16\biggl[\E[\yt_{t_1}\yt_{t_1+k_1}\ytt_{t_2}\ytt_{t_2+k_2}]-\E[\yt_{t_1},\yt_{t_1+k_1}]\E[\ytt_{t_2},\ytt_{t_2+k_2}]\biggr]\\
={}&16\biggl[\cum (\yt_{t_1},\yt_{t_1+k_1},\ytt_{t_2},\ytt_{t_2+k_2})+\frac{1}{144}\rho(t_2-t_1)\rho(t_2+k_2-(t_1+k_1))\\
&+\frac{1}{144}\rho(t_2+k_2-t_1)\rho(t_2-(t_1+k_1))\biggr],
}
where $(Y_t^{(j)})_{t\in\Z}:=\Big(I(X_t<X_t^{(j)})-\frac{1}{2}\Big)_{t\in\Z}$, $j=1,2$.

\item[(ii)] for Spearman's $\rho$
\als{\label{rho1}
&\Cov\biggl(\hunindr{k_1}{X_{t_1}}{X_{t_1+k_1}},\hunindr{k_2}{X_{t_2}}{X_{t_2+k_2}}\biggr)\notag\\
={}&4\sum_{\gamma\in\Gamma\{1,2,3\}}\sum_{\gt\in\Gamma\{1,4,5\}}\E\biggl[\biggl(I(X_{t_1}^{(\gamma(1))}<X_{t_1}^{(\gamma(2))})-\frac{1}{2}\biggr)\biggl(I(X_{t_1+k_1}^{(\gamma(1))}<X_{t_1+k_1}^{(\gamma(3))})-\frac{1}{2}\biggr)\notag\\
&\hspace{2cm} \biggl(I(X_{t_2}^{(\gt(1))}<X_{t_2}^{(\gt(2))})-\frac{1}{2}\biggr)\biggl(I(X_{t_2+k_2}^{(\gt(1))}<X_{t_2+k_2}^{(\gt(3))})-\frac{1}{2}\biggr)\biggr]\notag\\
={}&\sum_{\gamma\in\Gamma\{1,2,3\}}\sum_{\gt\in\Gamma\{1,4,5\}}4\,\cum\biggl(I(X_{t_1}^{(\gamma(1))}<X_{t_1}^{(\gamma(2))}),I(X_{t_1+k_1}^{(\gamma(1))}<X_{t_1+k_1}^{(\gamma(3))}),\notag\\
&\hspace{4cm}I(X_{t_2}^{(\gt(1))}<X_{t_2}^{(\gt(2))}),I(X_{t_2+k_2}^{(\gt(1))}<X_{t_2+k_2}^{(\gt(3))})\biggr)\notag\\
&+\frac{1}{9}\rho(t_2-t_1)\rho(t_2+k_2-(t_1+k_1))+\frac{1}{9}\rho(t_2+k_2-t_1)\rho(t_2-(t_1+k_1)),
}

where $\Gamma\{i,j,k\}$ denotes the set of all permutations of $\{i,j,k\}$.
\ei
\end{lem}
Lemma~\ref{auxlemmavar} is proved in Section~\ref{sec:proofauxlemmavar}. 


\subsubsection{Proof of (\ref{varblocktau})}
Let $\proct$ and $\proctt$ be independent copies of the strictly stationary $\beta$-mixing process $\proc$ that are independent of each other. Define processes $\procty$ and $\proctty$ by

\[
\procty=\Big(I(X_t<X_t^{(1)})-\frac{1}{2}\Big)_{t\in\Z}\qquad\text{ and }\qquad \proctty=\Big(I(X_t<X_t^{(2)})-\frac{1}{2}\Big)_{t\in\Z}.
\]

Note that the processes $\procty$ and $\proctty$ are strictly stationary.\\ 

By Lemma \ref{auxlemmavar} (i) we have for Kendall's $\tau$,
\begin{align*}
&\E\Big[\sum_{t_1=1}^{a_n}W_{n,t_1}^{\tau}(\omega)\sum_{t_2=1}^{a_n}W_{n,t_2}^{\tau}(\omega)\Big]\\
={}&\frac{1}{(2\pi)^2}\sum_{|k_1|\leq r_n}\sum_{|k_2|\leq r_n}\w{k_1}\w{k_2}e^{-i(k_1+k_2)\omega}\frac{4}{n^2}\sum_{t_1=1}^{a_n}\sum_{t_2=1}^{a_n}\E\Big[h_{1,k_1}^{\tau}(\bm{X}_{t_1,k_1})h_{1,k_2}^{\tau}(\bm{X}_{t_2,k_2})\Big]\\
={}&\frac{1}{(2\pi)^2}\sum_{|k_1|\leq r_n}\sum_{|k_2|\leq r_n}\w{k_1}\w{k_2}e^{-i(k_1+k_2)\omega}\frac{4}{n^2}\sum_{t_1=1}^{a_n}\sum_{t_2=1}^{a_n}16\Big[\cum(\yt_{t_1},\yt_{t_1+k_1},\ytt_{t_2},\ytt_{t_2+k_2})\\
&+\frac{1}{144}\rho(t_2-t_1)\rho(t_2+k_2-(t_1+k_1))+\frac{1}{144}\rho(t_2+k_2-t_1)\rho(t_2-(t_1+k_1))\Big].
\end{align*}
Next, let $\mathcal{T}_{k,a_n}:=\{t|t,t+k\in\{1,\dots,a_n\}\}$ and define for $j=1,2$,
\[
H_j^{\tau}(\omega):=\frac{1}{2\pi}\frac{2}{a_n}\sum_{|k|\leq r_n}\wk e^{-ik\omega}\sum_{t\in\mathcal{T}_{k,a_n}}4(Y^{(j)}_tY^{(j)}_{t+k}-\E[Y^{(j)}_{t}Y^{(j)}_{t+k}]).
\]
Then, similar arguments as in the proof of Lemma \ref{auxlemmavar} yield
\als{\label{covhtau}
\E[H_1^{\tau}(\omega)H_2^{\tau}(\omega)]={}&\frac{1}{(2\pi)^2}\sum_{|k_1|\leq r_n}\sum_{|k_2|\leq r_n}\w{k_1}\w{k_2}e^{-i(k_1+k_2)\omega}\frac{4}{n^2}\sum_{t_1\in\mathcal{T}_{k_1,a_n}}\sum_{t_2\in\mathcal{T}_{k_2,a_n}}\notag\\
&\quad16\Big[\cum(\yt_{t_1},\yt_{t_1+k_1},\ytt_{t_2},\ytt_{t_2+k_2})+\frac{1}{144}\rho(t_2-t_1)\rho(t_2+k_2-(t_1+k_1))\notag\\
&\qquad+\frac{1}{144}\rho(t_2+k_2-t_1)\rho(t_2-(t_1+k_1))\Big]
}
and consequently,
\al{
\Big|\E\Big[\sum_{t_1=1}^{a_n}W_{n,t_1}^{\tau}(\omega)\sum_{t_2=1}^{a_n}W_{n,t_2}^{\tau}(\omega)\Big]-\frac{a_n^2}{n^2}\E[H_1^{\tau}(\omega)H_2^{\tau}(\omega)]\Big|
\leq{}&|D_{1,n}|+|D_{2,n}|,
}
where 
\al{
D_{1,n}:={}&\frac{1}{(2\pi)^2}\sum_{|k_1|\leq r_n}\sum_{|k_2|\leq r_n}\w{k_1}\w{k_2}e^{-i(k_1+k_2)\omega}\frac{64}{n^2}\Big(\sum_{t_1=1}^{a_n}\sum_{t_2=1}^{a_n}-\sum_{t_1\in\mathcal{T}_{k_1,a_n}}\sum_{t_2\in\mathcal{T}_{k_2,a_n}}\Big)\\
&\cum(\yt_{t_1},\yt_{t_1+k_1},\ytt_{t_2},\ytt_{t_2+k_2}),
}

\al{
D_{2,n}:={}&\frac{1}{(2\pi)^2}\sum_{|k_1|\leq r_n}\sum_{|k_2|\leq r_n}\w{k_1}\w{k_2}e^{-i(k_1+k_2)\omega}\frac{64}{n^2}\Big(\sum_{t_1=1}^{a_n}\sum_{t_2=1}^{a_n}-\sum_{t_1\in\mathcal{T}_{k_1,a_n}}\sum_{t_2\in\mathcal{T}_{k_2,a_n}}\Big)\\
&\frac{1}{144}\Big\{\rho(t_2-t_1)\rho(t_2+k_2-(t_1+k_1))+\rho(t_2+k_2-t_1)\rho(t_2-(t_1+k_1))\Big\}.
}
We will now derive upper bounds for $|D_{1,n}|$ and $|D_{2,n}|$ separately. First, as $|w(\cdot)|\leq 1$, $|e^{i\cdot}|=1$ and $\mathcal{T}_{k,a_n}\subset\{1,\dots,a_n\}$, we have
\al{
|D_{1,n}|\leq{} &2\frac{1}{4\pi^2}\frac{64}{n^2}\sum_{|k_1|\leq r_n}\sum_{|k_2|\leq r_n}\sum_{t_1=1}^{a_n}\sum_{t_2=1}^{a_n}|\cum(Y_{t_1}^{(1)},Y_{t_1+k_1}^{(1)},Y_{t_2}^{(2)},Y_{t_2+k_2}^{(2)})|\\
\leq{} &\frac{32}{\pi^2n^2}\sum_{|u_3|\leq a_n+r_n}\sum_{|u_2|\leq a_n}\sum_{|u_1|\leq r_n}\sum_{t_1=1}^{a_n}|\cum(Y_{t_1}^{(1)},Y_{t_1+u_1}^{(1)},Y_{t_1+u_2}^{(2)},Y_{t_1+u_3}^{(2)})|\\
={} &\frac{32a_n}{\pi^2n^2}\sum_{|u_3|\leq a_n+r_n}\sum_{|u_2|\leq a_n}\sum_{|u_1|\leq r_n}|\cum(Y_{0}^{(1)},Y_{u_1}^{(1)},Y_{u_2}^{(2)},Y_{u_3}^{(2)})|,
}
where the latter inequality follows by the strict joint stationarity of the involved processes. Next, observe that, from Theorem 2.3.1 in \cite{brillinger75}, it follows that
\al{
\cum(Y_{t_1}^{(1)},Y_{t_2}^{(1)},Y_{t_3}^{(2)},Y_{t_4}^{(2)})
={}&\cum(I(X_{t_1}<X_{t_1}^{(1)}),I(X_{t_2}<X_{t_2}^{(1)}),I(X_{t_3}<X_{t_3}^{(2)}),I(X_{t_4}<X_{t_4}^{(2)}))\\
={}&\cum(I(V_{t_1}\in A_1),I(V_{t_2}\in A_2),I(V_{t_3}\in A_3),I(V_{t_4}\in A_4)),
}
where $V_{t_j}:=(X_{t_j},X_{t_j}^{(1)},X_{t_j}^{(2)})$, $A_1=\{x\in\R^3:x_1<x_2\}=A_2$ and $A_3=\{x\in\R^3:x_1<x_3\}=A_4$. Furthermore, let $u_0:=0$ and consider the set
\[
\mathcal{S}_m:=\Big\{(u_1,\dots,u_p)\in\Z^p|\max_{i,j=0,\dots,p}|u_i-u_j|=m\Big\}
\]
whose cardinality is $\leq c_p(m+1)^{p-1}$. Hence, applying Lemma \ref{cum} with $q=2$ yields 
\al{
&\sum_{|u_3|\leq a_n+r_n}\sum_{|u_2|\leq a_n}\sum_{|u_1|\leq r_n}|\cum(Y_{0}^{(1)},Y_{u_1}^{(1)},Y_{u_2}^{(2)},Y_{u_3}^{(2)})|\\
\leq{} &C_{4,2}\sum_{|u_3|\leq a_n+r_n}\sum_{|u_2|\leq a_n}\sum_{|u_1|\leq r_n}\alpha(\max_{i,j=0,1,2,3}|u_i-u_j|)\\
\leq{} &C_{4,2}\sum_{m=0}^{\infty}\sum_{(u_1,u_2,u_3)\in\mathcal{S}_m}\alpha(m)
\leq{} C_{4,2}c_3(1+\sum_{m=1}^{\infty}m^2\alpha(m))
={} O(1),
}
where we used Assumption (N1) for the last estimate and $\alpha:=\alpha^{X}$ are the mixing coefficients of the process $\proc$. Therefore,
\als{\label{D1n}
D_{1,n}=O\Big(\frac{a_n}{n^2}\Big)=o\Big(\frac{r_n}{a_n}\Big).
}
Next, 
\al{
D_{2,n}={}&\frac{1}{(2\pi)^2}\sum_{|k_1|\leq r_n}\sum_{|k_2|\leq r_n}\w{k_1}\w{k_2}e^{-i(k_1+k_2)\omega}\frac{64}{n^2}\Big(\sum_{t_1=1}^{a_n}-\sum_{t_1\in\mathcal{T}_{k_1,a_n}}\Big)\sum_{t_2=1}^{a_n}\\
&\hspace{1cm}\frac{1}{144}\Big\{\rho(t_2-t_1)\rho(t_2+k_2-(t_1+k_1))+\rho(t_2+k_2-t_1)\rho(t_2-(t_1+k_1))\Big\}\\
&+\frac{1}{(2\pi)^2}\sum_{|k_1|\leq r_n}\sum_{|k_2|\leq r_n}\w{k_1}\w{k_2}e^{-i(k_1+k_2)\omega}\frac{64}{n^2}\Big(\sum_{t_2=1}^{a_n}-\sum_{t_2\in\mathcal{T}_{k_2,a_n}}\Big)\sum_{t_1\in\mathcal{T}_{k_1,a_n}}\\
&\hspace{1cm}\frac{1}{144}\Big\{\rho(t_2-t_1)\rho(t_2+k_2-(t_1+k_1))+\rho(t_2+k_2-t_1)\rho(t_2-(t_1+k_1))\Big\}\\
={}&D_{2,n}^{(1)}+D_{2,n}^{(2)}
}
and observing that $(\sum_{t_1=1}^{a_n}-\sum_{t_1\in\mathcal{T}_{k_1,a_n}})$ contains $O(r_n)$ summands, we obtain by assumption (N2),
\al{
|D_{2,n}^{(1)}|\leq{}&\frac{1}{(2\pi)^2}\frac{1}{144}\frac{64}{n^2}\sum_{|k_1|\leq r_n}\Big(\sum_{t_1=1}^{a_n}-\sum_{t_1\in\mathcal{T}_{k_1,a_n}}\Big)\Big[\sum_{t_2=1}^{a_n}|\rho(t_2-t_1)|\sum_{|k_2|\leq r_n}|\rho(t_2+k_2-(t_1+k_1))|\\
&+\sum_{t_2=1}^{a_n}|\rho(t_2-(t_1+k_1))|\sum_{|k_2|\leq r_n}|\rho(t_2+k_2-t_1)|\Big]\\
\leq{}&2\frac{1}{(2\pi)^2}\frac{1}{144}\frac{64}{n^2}r_n^2\sum_{t_2\in\Z}|\rho(t_2)|\sum_{|k_2|\leq r_n}|\rho(k_2)
={}O\Big(\frac{r_n^2}{n^2}\Big).
}
Analogously, $D_{2,n}^{(2)}=O\Big(\frac{r_n^2}{n^2}\Big)$ and hence,
\als{\label{D2n}
D_{2,n}=O\Big(\frac{r_n^2}{n^2}\Big)=o\Big(\frac{r_n}{a_n}\Big).
}
Then, equations (\ref{D1n}) and (\ref{D2n}) together yield
\als{\label{diffWh}
\Big|\E\Big[\sum_{t_1=1}^{a_n}W_{n,t_1}^{\tau}(\omega)\sum_{t_2=1}^{a_n}W_{n,t_2}^{\tau}(\omega)\Big]-\frac{a_n^2}{n^2}\E[H_1^{\tau}(\omega)H_2^{\tau}(\omega)]\Big|=o\Big(\frac{r_n}{a_n}\Big).
}
Next, observe that $h_j^{\tau}(\omega)$, $j=1,2$, is eight times the classical centered lag-window estimator of the spectral density of the stationary process $(Y_t^{(j)})_{t\in\Z}$ based on the observations $Y_1^{(j)},\dots,Y_{a_n}^{(j)}$. Consequently, if we show that 
\als{\label{boundcum}
&\frac{16}{\pi^2 a_n^2}\sum_{|k_1|\leq r_n}\sum_{|k_2|\leq r_n}\w{k_1}\w{k_2}e^{-i(k_1+k_2)\omega}\sum_{t_1\in\mathcal{T}_{k_1,a_n}}\sum_{t_2\in\mathcal{T}_{k_2,a_n}}|\cum(Y_{t_1}^{(1)},Y_{t_1+k_1}^{(1)},Y_{t_2}^{(2)},Y_{t_2+k_2}^{(2)})|
=o\Big(\frac{r_n}{a_n}\Big)
} 
the same arguments as given in the proof of Theorem 9.3.4 in \cite{anderson71} yield
\als{\label{varhtau}
\E[H_1^{\tau}(\omega)H_2^{\tau}(\omega)]=\frac{4}{9}\frac{r_n}{a_n}\spec_{\rho}^2(\omega)\int_{-1}^1w^2(u)du(1+I(\omega\in\Big\{0,\pm\pi\Big\}))+o\Big(\frac{r_n}{a_n}\Big).
}

Finally, (\ref{boundcum}) can be proved using similar arguments and equation (\ref{varhtau}) together with equation (\ref{diffWh}) conclude the proof of (\ref{varblocktau}).\\

\subsubsection{Proof of (\ref{varblockrho})}
Let $(X_t^{(j)})_{t\in\Z}$, $j=1,\dots,5$ be independent copies of the strictly stationary process $\proc$ that are independent of each other. Then, by Lemma \ref{auxlemmavar} (ii) for Spearman's $\rho$,
\als{\label{spearmandec}
&\E\Big[\sum_{t_1=1}^{a_n}W_{n,t_1}^{\rho}(\omega)\sum_{t_2=1}^{a_n}W_{n,t_2}^{\rho}(\omega)\Big]\notag\\
={}&\frac{1}{(2\pi)^2}\sum_{|k_1|\leq r_n}\sum_{|k_2|\leq r_n}\w{k_1}\w{k_2}e^{-i(k_1+k_2)\omega}\frac{9}{n^2}\sum_{t_1=1}^{a_n}\sum_{t_2=1}^{a_n}\E\Big[h_{1,k_1}^{\rho}(\bm{X}_{t_1,k_1})h_{1,k_2}^{\rho}(\bm{X}_{t_2,k_2})\Big]\notag\\
={}&\frac{1}{(2\pi)^2}\sum_{|k_1|\leq r_n}\sum_{|k_2|\leq r_n}\w{k_1}\w{k_2}e^{-i(k_1+k_2)\omega}\frac{9}{n^2}\sum_{t_1=1}^{a_n}\sum_{t_2=1}^{a_n}\Big[\Big\{
4\sum_{\gamma\in\Gamma\{1,2,3\}}\sum_{\gt\in\Gamma\{1,4,5\}}\notag\\
&\cum\Big(I(X_{t_1}^{(\gamma(1))}<X_{t_1}^{(\gamma(2))}),I(X_{t_1+k_1}^{(\gamma(1))}<X_{t_1+k_1}^{(\gamma(3))}),I(X_{t_2}^{(\gt(1))}<X_{t_2}^{(\gt(2))}),I(X_{t_2+k_2}^{(\gt(1))}<X_{t_2+k_2}^{(\gt(3))})\Big)\Big\}\notag\\
&+\frac{1}{9}\rho(t_2-t_1)\rho(t_2+k_2-(t_1+k_1))+\frac{1}{9}\rho(t_2+k_2-t_1)\rho(t_2-(t_1+k_1))\Big],
}
where $\Gamma\{i,j,k\}$ denotes the set of all permutations of $\{i,j,k\}$.\\

Next, let $\mathcal{T}_{k,a_n}:=\{t|t,t+k\in\{1,\dots,a_n\}\}$ and define
\al{
H_1^{\rho}(\omega):=&\frac{1}{2\pi}\frac{3}{a_n}\sum_{|k|\leq r_n}\wk e^{-ik\omega}\sum_{t\in\mathcal{T}_{k,a_n}}\Big[\sum_{\gamma\in\Gamma\{1,2,3\}}\\
&\hspace{2cm}2\Big(I(X_{t}^{(\gamma(1))}<X_{t}^{(\gamma(2))})-\frac{1}{2}\Big)\Big(I(X_{t+k}^{(\gamma(1))}<X_{t+k}^{(\gamma(3))})-\frac{1}{2}\Big)-\rho(k)\Big]
}
and
\al{
H_2^{\rho}(\omega):=&\frac{1}{2\pi}\frac{3}{a_n}\sum_{|k|\leq r_n}\wk e^{-ik\omega}\sum_{t\in\mathcal{T}_{k,a_n}}\Big[\sum_{\gt\in\Gamma\Big\{1,4,5\Big\}}\\
&\hspace{2cm}2\Big(I(X_{t}^{(\gt(1))}<X_{t}^{(\gt(2))})-\frac{1}{2}\Big)\Big(I(X_{t+k}^{(\gt(1))}<X_{t+k}^{(\gt(3))})-\frac{1}{2}\Big)-\rho(k)\Big].
}
Then, similarly as in the proof of Lemma \ref{auxlemmavar},
\als{\label{covhrho}
&\E[H_1^{\rho}(\omega)H_2^{\rho}(\omega)]\notag\\
={}&\frac{1}{(2\pi)^2}\sum_{|k_1|\leq r_n}\sum_{|k_2|\leq r_n}\w{k_1}\w{k_2}e^{-i(k_1+k_2)\omega}\frac{1}{a_n^2}\sum_{t_1\in\mathcal{T}_{k_1,a_n}}\sum_{t_2\in\mathcal{T}_{k_2,a_n}}\Big[\Big\{
36\sum_{\gamma\in\Gamma\{1,2,3\}}\sum_{\gt\in\Gamma\{1,4,5\}}\notag\\
&\cum\Big(I(X_{t_1}^{(\gamma(1))}<X_{t_1}^{(\gamma(2))}),I(X_{t_1+k_1}^{(\gamma(1))}<X_{t_1+k_1}^{(\gamma(3))}),I(X_{t_2}^{(\gt(1))}<X_{t_2}^{(\gt(2))}),I(X_{t_2+k_2}^{(\gt(1))}<X_{t_2+k_2}^{(\gt(3))})\Big)\Big\}\notag\\
&+\rho(t_2-t_1)\rho(t_2+k_2-(t_1+k_1))+\rho(t_2+k_2-t_1)\rho(t_2-(t_1+k_1))\Big]
}
and analogous arguments as for Kendall's $\tau$ give
\als{\label{diffWhrho}
\Big|\E\Big[\sum_{t_1=1}^{a_n}W_{n,t_1}^{\rho}(\omega)\sum_{t_2=1}^{a_n}W_{n,t_2}^{\rho}(\omega)\Big]-\frac{a_n^2}{n^2}\E[H_1^{\rho}(\omega)H_2^{\rho}(\omega)]\Big|=o\Big(\frac{r_n}{a_n}\Big)
}
Next, similar arguments as were used in order to derive (\ref{D1n}) yield  
\al{
&\Big|\frac{1}{(2\pi)^2}\frac{36}{a_n^2}\sum_{|k_1|\leq r_n}\sum_{|k_2|\leq r_n}\w{k_1}\w{k_2}e^{-i(k_1+k_2)\omega}\sum_{t_1\in\mathcal{T}_{k_1,a_n}}\sum_{t_2\in\mathcal{T}_{k_2,a_n}}\sum_{\gamma\in\Gamma\{1,2,3\}}\sum_{\gt\in\Gamma\{1,4,5\}}\notag\\
&\cum\Big(I(X_{t_1}^{(\gamma(1))}<X_{t_1}^{(\gamma(2))}),I(X_{t_1+k_1}^{(\gamma(1))}<X_{t_1+k_1}^{(\gamma(3))}),I(X_{t_2}^{(\gt(1))}<X_{t_2}^{(\gt(2))}),I(X_{t_2+k_2}^{(\gt(1))}<X_{t_2+k_2}^{(\gt(3))})\Big)\Big|\\
&\hspace{12cm}=o\Big(\frac{r_n}{a_n}\Big),
}
and the same arguments as in the proof of Theorem 9.3.4 in \cite{anderson71} show
\als{\label{varhrho}
\E[H_1^{\rho}(\omega)H_2^{\rho}(\omega)]=\frac{r_n}{a_n}\spec_{\rho}^2(\omega)\int_{-1}^1w^2(u)du(1+I(\omega\in\Big\{0,\pm\pi\Big\}))+o\Big(\frac{r_n}{a_n}\Big).
}
Hence, equation (\ref{varhrho}) together with equation (\ref{diffWhrho}) conclude the proof of (\ref{varblockrho}).

\subsubsection{Proof of (\ref{fourthmoment}) and (\ref{lowerboundvar})}
By Lemma \ref{auxlemmavar}, we know that for Kendall's tau
\bea
\E\Big[\prod_{j=1}^4h_{1,k_j}^{\tau}\colvec{X_{t_j}}{X_{t_j+k_j}}\Big]&=&\E\Big[\prod_{j=1}^4\Big(h^{\tau}\Big(\colvec{X_{t_j}}{X_{t_j+k_j}},\colvec{X_{t_j}^{(j)}}{X_{t_j+k_j}^{(j)}}\Big)-\tau_{k_j}\Big)\Big]\\
&=&\E\Big[\prod_{j=1}^4(4Y_{t_j}^{(j)}Y_{t_j+k_j}^{(j)}-\tau_{k_j})\Big],
\eea
where $(Y_{t}^{(j)})_{t\in\Z}=\Big(I(X_t<X_t^{(j)})-\frac{1}{2}\Big)_{t\in\Z}$. Therefore, we can write
\als{\label{zeta=xi}
\E\Big[\Big(\sum_{t\in\Gamma_j}\zeta^{\tau}_{n,t}(\omega)\Big)^4\Big]={}&\E\Big[\sum_{t_1,t_2,t_3,t_4\in\Gamma_j}W^{\tau}_{n,t_1}(\omega)W^{\tau}_{n,t_2}(\omega)W^{\tau}_{n,t_3}(\omega)W^{\tau}_{n,t_4}(\omega)\Big]\notag\\
={}&\E\Big[\sum_{t_1,t_2,t_3,t_4\in\Gamma_j}\vartheta^{\tau}_{n,t_1}(\omega)\vartheta^{\tau}_{n,t_2}(\omega)\vartheta^{\tau}_{n,t_3}(\omega)\vartheta^{\tau}_{n,t_4}(\omega)\Big],
}
where $\vartheta^{\tau}_{t_l}=\frac{1}{2\pi}\sum_{|k_l|\leq r_n}w\Big(\frac{k_l}{r_n}\Big)e^{-ik_l\omega}\frac{2}{n}[4Y_{t_l}^{(l)}Y_{t_l+k_l}^{(l)}-\tau_{k_l}]$, $l=1,\dots,4$. Observing that by construction $\E[\sum_{t_l\in\Gamma_j}\vartheta^{\tau}_{t_l}(\omega)]=0$, we express the fourth moment in terms of a fourth order cumulant and 3 products of second order cumulants, that is 
\al{
\E\Big[\sum_{t_1,t_2,t_3,t_4\in\Gamma_j}&\vartheta^{\tau}_{n,t_1}(\omega)\vartheta^{\tau}_{n,t_2}(\omega)\vartheta^{\tau}_{n,t_3}(\omega)\vartheta^{\tau}_{n,t_4}(\omega)\Big]=\cum\Big(\sum_{t_l\in\Gamma_j}\vartheta^{\tau}_{n,t_l}(\omega);l=1,\dots,4\Big)\\
&+\cum\Big(\sum_{t_l\in\Gamma_j}\vartheta^{\tau}_{n,t_l}(\omega);l=1,2\Big)\cum\Big(\sum_{t_l\in\Gamma_j}\vartheta^{\tau}_{n,t_l}(\omega);l=3,4\Big)\\
&\quad+\cum\Big(\sum_{t_l\in\Gamma_j}\vartheta^{\tau}_{n,t_l}(\omega);l=1,3\Big)\cum\Big(\sum_{t_l\in\Gamma_j}\vartheta^{\tau}_{n,t_l}(\omega);l=2,4\Big)\\
&\qquad+\cum\Big(\sum_{t_l\in\Gamma_j}\vartheta^{\tau}_{n,t_l}(\omega);l=1,4\Big)\cum\Big(\sum_{t_l\in\Gamma_j}\vartheta^{\tau}_{n,t_l}(\omega);l=2,3\Big)
}
Note that by construction for all $k,l\in\{1,\dots,4\}$, $\sum_{t_k\in\Gamma_j}\vartheta^{\tau}_{n,t_k}\overset{\mathcal{D}}{=}\sum_{t_l\in\Gamma_j}\vartheta^{\tau}_{n,t_l}$. Therefore, each of the second order cumulants is equal to
\al{
\cum\Big(\sum_{t_l\in\Gamma_j}\vartheta^{\tau}_{n,t_l}(\omega);l=1,2\Big)={}&\E\Big[\sum_{t_1\in\Gamma_j}\vartheta^{\tau}_{n,t_1}(\omega)\sum_{t_2\in\Gamma_j}\vartheta^{\tau}_{n,t_2}(\omega)\Big]\\
={}&\E\Big[\Big(\sum_{t\in\Gamma_j}W^{\tau}_{n,t}(\omega)\Big)^2\Big]=\Big(\Var\Big(\sum_{t\in\Gamma_j}W^{\tau}_{n,t}(\omega)\Big)\Big)^2,
}
where we have used a similar argument as in equation (\ref{zeta=xi}). Hence, we obtain by Theorem 2.3.1 in \cite{brillinger75} 
\al{
&\E\Big[\Big(\sum_{t_l\in\Gamma_j}\vartheta^{\tau}_{n,t_l}(\omega)\Big)^4\Big]
={}\sum_{t_1,t_2,t_3,t_4\in\Gamma_j}\cum(\vartheta^{\tau}_{n,t_1}(\omega),\vartheta^{\tau}_{n,t_2}(\omega),\vartheta^{\tau}_{n,t_3}(\omega),\vartheta^{\tau}_{n,t_4}(\omega))+3\Big(\Var\Big(\sum_{t\in\Gamma_j}W^{\tau}_{n,t}(\omega)\Big)\Big)^2
}
Following the arguments of \cite{rosenblatt84} on pages 1177-1178, we can express the fourth order cumulant of products in terms of cumulants of the factors, i.e. we obtain
\al{
&\sum_{t_1,t_2,t_3,t_4\in\Gamma_j}\cum(\vartheta^{\tau}_{n,t_1}(\omega),\vartheta^{\tau}_{n,t_2}(\omega),\vartheta^{\tau}_{n,t_3}(\omega),\vartheta^{\tau}_{n,t_4}(\omega))\\
={}&\frac{1}{(2\pi)^4}\frac{2^44^4}{n^4}\sum_{|k_1|,|k_2|,|k_3|,|k_4|\leq r_n}\sum_{t_1,t_2,t_3,t_4\in\Gamma_j}\Big(\prod_{l=1}^4 w\Big(\frac{k_l}{r_n}\Big)e^{-ik_l\omega}\Big)\cum(Y_{t_l}^{(l)}Y_{t_l+k_l}^{(l)};l=1,\dots,4)\\
={}&\frac{1}{(2\pi)^4}\frac{2^44^4}{n^4}\sum_{|k_1|,|k_2|,|k_3|,|k_4|\leq r_n}\sum_{t_1,t_2,t_3,t_4\in\Gamma_j}\Big(\prod_{l=1}^4 w\Big(\frac{k_l}{r_n}\Big)e^{-ik_l\omega}\Big)\\
&\qquad\qquad\qquad\qquad\qquad\cdot\sum_{\nu}\cum(Y_s^{(l)},s\in\nu_1)\cdots\cum(Y_s^{(l)},s\in\nu_r)
}
where the latter sum extends over all indecomposable partitions $\nu=\nu_1\cup\dots\cup\nu_r$ of the table
\begin{align*}
&t_1\qquad t_1+k_1\\
&t_2\qquad t_2+k_2\\
&t_3\qquad t_3+k_3\\
&t_4\qquad t_4+k_4.
\end{align*}
In order to bound this sum, we need that for $2\leq p\leq 8$ and $t\in\Z$
\als{\label{summabilitycum8}
\sum_{u_1,\dots,u_{p-1}\in\Z}|\cum(I(V_{t}\in A_1),I(V_{t+u_1}\in A_2),\dots,I(V_{t+u_{p-1}}\in A_{p}))|<\infty,
}
with $V_t:=(X_t,X_t^{(1)},\dots,X_t^{(4)})$ and measurable sets $A_1,\dots,A_p\subset\R^5$. This follows by  Lemma \ref{cum} and Assumption (N1): 
\al{
&\sum_{u_1,\dots,u_{p-1}\in\Z}|\cum(I(V_{t}\in A_1),I(V_{t+u_1}\in A_2),\dots,I(V_{t+u_{p-1}}\in A_{p}))|\\
\leq{} &C_{p,4}\sum_{u_1,\dots,u_{p-1}\in\Z}\alpha(\max_{i,j=0,\dots,p-1}|u_i-u_j|)\\
\leq{} &C_{p,4}\sum_{m=0}^{\infty}\sum_{u_1,\dots,u_{p-1}\in\mathcal{S}_m}\alpha(m)
\leq{} C_{p,4}c_{p-1}(1+\sum_{m=1}^{\infty}m^{p-2}\alpha(m))<\infty,
}
where $u_0:=0$ and $\mathcal{S}_m:=\{(u_1,\dots,u_p)\in\Z^p|\max_{i,j=0,\dots,p}|u_i-u_j|=m\}$ have been introduced in the proof of (\ref{varblocktau}).\\
 
Next, arguments as in \cite{rosenblatt84} on page 1177--1178 yield
\bea
\sum_{t_1,t_2,t_3,t_4\in\Gamma_j}|\cum(\vartheta^{\tau}_{n,t_1}(\omega),\vartheta^{\tau}_{n,t_2}(\omega),\vartheta^{\tau}_{n,t_3}(\omega),\vartheta^{\tau}_{n,t_4}(\omega))|=O\Big(\frac{p_n^2r_n^2}{n^4}\Big)
\eea
and together with (\ref{varblocktau}) we obtain
\bea
\sum_{j=1}^{\mu_n}\E\Big[\Big(\sum_{t\in\Gamma_j}\zeta^{\tau}_{n,t}(\omega)\Big)^4\Big]=O\Big(\frac{\mu_np_n^2r_n^2}{n^4}\Big).
\eea
Furthermore, as $\spec_{\xi}(\omega)\neq 0$, by (\ref{varblocktau}) 
\bea
\Big(\sum_{j=1}^{\mu_n}\Var\Big(\sum_{t\in\Gamma_j}\zeta^{\tau}_{n,t}(\omega)\Big)\Big)^2\geq c\frac{\mu_n^2p_n^2r_n^2}{n^4}
\eea
for some constant $c>0$ and $n$ sufficiently large. This yields (\ref{fourthmoment}) and (\ref{lowerboundvar}) in the case where $\xi$ is Kendall's $\tau$.\\

In the case where $\xi$ is Spearman's $\rho$ we have by Lemma \ref{auxlemmavar} 
\als{\label{zeta=xirho}
\E\Big[\Big(\sum_{t\in\Gamma_j}\zeta^{\rho}_{n,t}(\omega)\Big)^4\Big]={}&\E\Big[\sum_{t_1,t_2,t_3,t_4\in\Gamma_j}W^{\rho}_{n,t_1}(\omega)W^{\rho}_{n,t_2}(\omega)W^{\rho}_{n,t_3}(\omega)W^{\rho}_{n,t_4}(\omega)\Big]\notag\\
={}&\E\Big[\sum_{t_1,t_2,t_3,t_4\in\Gamma_j}\vartheta^{\rho}_{n,t_1}(\omega)\vartheta^{\rho}_{n,t_2}(\omega)\vartheta^{\rho}_{n,t_3}(\omega)\vartheta^{\rho}_{n,t_4}(\omega)\Big],
}
where 
\al{
\vartheta^{\rho}_{t_l}={}&\frac{1}{2\pi}\sum_{|k_l|\leq r_n}w\Big(\frac{k_l}{r_n}\Big)e^{-ik_l\omega}\frac{3}{n}\Big[\sum_{\gamma_l\in\Gamma\{1,2l,2l+1\}}2\Big(I(X_{t_l}^{(\gamma_l(1))}<X_{t_l}^{(\gamma_l(2))})-\frac{1}{2}\Big)\\
&\hspace{4cm}\cdot\Big(I(X_{t_l+k_l}^{(\gamma_l(1))}<X_{t_l+k_l}^{(\gamma_l(3))})-\frac{1}{2}\Big)-\rho_{k_l}\Big],\qquad l=1,\dots,4.
}
Observing that by construction $\E[\sum_{t_l\in\Gamma_j}\vartheta^{\rho}_{t_l}(\omega)]=0$, we have similarly as for Kendall's $\tau$ 
\al{
&\E\Big[\Big(\sum_{t_l\in\Gamma_j}\vartheta^{\rho}_{n,t_l}(\omega)\Big)^4\Big]\\
={}&\sum_{t_1,t_2,t_3,t_4\in\Gamma_j}\cum(\vartheta^{\rho}_{n,t_1}(\omega),\vartheta^{\rho}_{n,t_2}(\omega),\vartheta^{\rho}_{n,t_3}(\omega),\vartheta^{\rho}_{n,t_4}(\omega))+3\Big(\Var\Big(\sum_{t\in\Gamma_j}W^{\rho}_{n,t}(\omega)\Big)\Big)^2,
}
where, 
\als{\label{zerlegungcumprod}
&\sum_{t_1,t_2,t_3,t_4\in\Gamma_j}\cum_4(\vartheta^{\rho}_{n,t_1}(\omega),\vartheta^{\rho}_{n,t_2}(\omega),\vartheta^{\rho}_{n,t_3}(\omega),\vartheta^{\rho}_{n,t_4}(\omega))\notag\\
={}&\frac{1}{(2\pi)^4}\frac{3^42^4}{n^4}\sum_{|k_1|,|k_2|,|k_3|,|k_4|\leq r_n}\sum_{t_1,t_2,t_3,t_4\in\Gamma_j}\Big(\prod_{l=1}^4 w\Big(\frac{k_l}{r_n}\Big)e^{-ik_l\omega}\sum_{\gamma_l\in\Gamma\{1,2l,2l+1\}}\Big)\notag\\
&\qquad\cum_4\Big(I(X_{t_l}^{(\gamma_l(1))}<X_{t_l}^{(\gamma_l(2))})I(X_{t_l+k_l}^{(\gamma_l(1))}<X_{t_l+k_l}^{(\gamma_l(3))});l=1,\dots,4\Big).
}
Following the arguments of \cite{rosenblatt84} on pages 1177-1178, we express the fourth order cumulants of products of random variables in terms of cumulants of the factors which, similarly as in (\ref{summabilitycum8}), can be bounded by Lemma \ref{cum} for $V_t:=(X_t^{(1)},\dots,X_t^{(9)})$ and $A_1,\dots,A_p\in\R^{9}$. After that, arguments as in \cite{rosenblatt84} yield
\bea
\sum_{t_1,t_2,t_3,t_4\in\Gamma_j}|\cum(\vartheta^{\rho}_{n,t_1}(\omega),\vartheta^{\rho}_{n,t_2}(\omega),\vartheta^{\rho}_{n,t_3}(\omega),\vartheta^{\rho}_{n,t_4}(\omega))|=O\Big(\frac{p_n^2r_n^2}{n^4}\Big)
\eea
and together with (\ref{varblockrho}) we obtain
$
\sum_{j=1}^{\mu_n}\E\big[\big(\sum_{t\in\Gamma_j}\zeta^{\rho}_{n,t}(\omega)\big)^4\big]=O\big(\frac{\mu_np_n^2r_n^2}{n^4}\big).
$
Furthermore, as $\spec_{\xi}(\omega)\neq 0$, by (\ref{varblockrho}),
\bea
\Big(\sum_{j=1}^{\mu_n}\Var\Big(\sum_{t\in\Gamma_j}\zeta^{\rho}_{n,t}(\omega)\Big)\Big)^2\geq c\frac{\mu_n^2p_n^2r_n^2}{n^4}
\eea
for some constant $c>0$ and $n$ sufficiently large. This concludes the proof. \hfill $\Box$

\bigskip

{\bf Acknowledgements.} The authors would like to thank Martina Stein, who typed parts of this manuscript with considerable technical expertise.
This work has been supported in part by the Collaborative
Research Center ``Statistical modeling of nonlinear dynamic processes'' (SFB 823, Teilprojekt A1, C1) of the German Research Foundation (DFG).

 \setlength{\bibsep}{3pt}
\bibliographystyle{apalike}
\bibliography{Uspektraldichte}

\begin{thebibliography}{}

\bibitem[Ahdesm{\"a}ki et~al., 2005]{ahdesmaki2005}
Ahdesm{\"a}ki, M., L{\"a}hdesm{\"a}ki, H., Pearson, R., Huttunen, H., and
  Yli-Harja, O. (2005).
\newblock Robust detection of periodic time series measured from biological
  systems.
\newblock {\em BMC bioinformatics}, 6(1):1.

\bibitem[Anderson, 1971]{anderson71}
Anderson, T.~W. (1971).
\newblock {\em The Statistical Analysis of Time Series}.
\newblock John Wiley and Sons, Inc.

\bibitem[Arcones and Yu, 1994]{arcones94}
Arcones, M.~A. and Yu, B. (1994).
\newblock Central limit theorems for empirical and u-processes of stationary
  mixing sequences.
\newblock {\em Journal of Theoretical Probability}, 7(1):47--71.

\bibitem[Berbee, 1979]{berbee79}
Berbee, H.~C. (1979).
\newblock Random walks with stationary increments and renewal theory.
\newblock {\em MC Tracts}, 112:1--223.

\bibitem[Birr et~al., 2014]{birr2014}
Birr, S., Volgushev, S., Kley, T., Dette, H., and Hallin, M. (2014).
\newblock Quantile spectral analysis for locally stationary time series.
\newblock {\em arXiv preprint arXiv:1404.4605}.

\bibitem[Blomqvist, 1950]{blomqvist1950}
Blomqvist, N. (1950).
\newblock On a measure of dependence between two random variables.
\newblock {\em The Annals of Mathematical Statistics}, 21(4):593--600.

\bibitem[Bradley, 2005]{bradley05}
Bradley, R.~C. (2005).
\newblock Basic properties of strong mixing conditions. a survey and some open
  questions.
\newblock {\em Probability Surveys}, 2:107--144.

\bibitem[Brillinger, 1975]{brillinger75}
Brillinger, D.~R. (1975).
\newblock {\em Time series: Data Analysis and Theory}.
\newblock Holt, Rinehart and Winston, Inc.

\bibitem[Brockwell and Davis, 1987]{BrockwellDavis1987}
Brockwell, P.~J. and Davis, R.~A. (1987).
\newblock {\em Time Series: Theory and Methods}.
\newblock Springer Series in Statistics. Springer, New York.

\bibitem[Carcea and Serfling, 2015]{carcea2015}
Carcea, M. and Serfling, R. (2015).
\newblock A gini autocovariance function for time series modelling.
\newblock {\em Journal of Time Series Analysis}, 36(6):817--838.

\bibitem[Davis et~al., 2013]{DavisMikoschZhao2013}
Davis, R.~A., Mikosch, T., and Zhao, Y. (2013).
\newblock Measures of serial extremal dependence and their estimation.
\newblock {\em Stochastic Processes and their Applications}, 123:2575--2602.

\bibitem[Dette et~al., 2015]{DetteEtAl2013}
Dette, H., Hallin, M., Kley, T., Volgushev, S., et~al. (2015).
\newblock Of copulas, quantiles, ranks and spectra: An $l_1$-approach to
  spectral analysis.
\newblock {\em Bernoulli}, 21(2):781--831.

\bibitem[Eberlein, 1984]{eberlein84}
Eberlein, E. (1984).
\newblock Weak convergence of partial sums of absolutely regular sequences.
\newblock {\em Statistics \& Probability Letters}, 2(5):291--293.

\bibitem[Hagemann, 2013]{Hagemann2013}
Hagemann, A. (2013).
\newblock Robust spectral analysis (arxiv:1111.1965v2).
\newblock {\em ArXiv e-prints}.

\bibitem[Hong, 1999]{hong1999}
Hong, Y. (1999).
\newblock Hypothesis testing in time series via the empirical characteristic
  function: a generalized spectral density approach.
\newblock {\em Journal of the American Statistical Association},
  94(448):1201--1220.

\bibitem[Hong, 2000]{hong2000}
Hong, Y. (2000).
\newblock Generalized spectral tests for serial dependence.
\newblock {\em Journal of the Royal Statistical Society: Series B (Statistical
  Methodology)}, 62(3):557--574.

\bibitem[Kallenberg, 2010]{kallenberg10}
Kallenberg, O. (2010).
\newblock {\em Foundations of modern probability}.
\newblock Springer New York, 2 edition.

\bibitem[Kley, 2014]{kley14}
Kley, T. (2014).
\newblock {\em Quantile-Based Spectral Analysis: Asymptotic Theory and
  Computation}.
\newblock PhD thesis, Ruhr-Universit\"at Bochum.

\bibitem[Kley et~al., 2016]{DetteEtAl2016}
Kley, T., Volgushev, S., Dette, H., Hallin, M., et~al. (2016).
\newblock Quantile spectral processes: Asymptotic analysis and inference.
\newblock {\em Bernoulli}, 22(3):1770--1807.

\bibitem[Lee, 1990]{lee90}
Lee, A. (1990).
\newblock U-statistics, volume 110 of statistics: Textbooks and monographs.

\bibitem[Li, 2008]{Li2008}
Li, T.-H. (2008).
\newblock Laplace periodogram for time series analysis.
\newblock {\em Journal of the American Statistical Association},
  103(482):757--768.

\bibitem[Li, 2012]{Li2012}
Li, T.-H. (2012).
\newblock Quantile periodograms.
\newblock {\em Journal of the American Statistical Association},
  107(498):765--776.

\bibitem[Rosenblatt, 1984]{rosenblatt84}
Rosenblatt, M. (1984).
\newblock Asymptotic normality, strong mixing and spectral density estimates.
\newblock {\em The Annals of Probability}, pages 1167--1180.

\bibitem[Schechtman and Yitzhaki, 1987]{Schechtman1987}
Schechtman, E. and Yitzhaki, S. (1987).
\newblock A measure of association based on gini's mean difference.
\newblock {\em Communications in Statistics - Theory and Methods},
  16(1):207--231.

\bibitem[Schmid et~al., 2010]{schmid10}
Schmid, F., Schmidt, R., Blumentritt, T., Gai{\ss}er, S., and Ruppert, M.
  (2010).
\newblock Copula-based measures of multivariate association.
\newblock In {\em Copula theory and its applications}, pages 209--236.
  Springer.

\bibitem[Wald and Wolfowitz, 1943]{WaldWolfowitz1943}
Wald, A. and Wolfowitz, J. (1943).
\newblock An exact test for randomness in the non-parametric case based on
  serial correlation.
\newblock {\em The Annals of Mathematical Statistics}, 14(4):378--388.

\bibitem[Yoshihara, 1976]{yoshihara76}
Yoshihara, K.-i. (1976).
\newblock Limiting behavior of u-statistics for stationary, absolutely regular
  processes.
\newblock {\em Probability Theory and Related Fields}, 35(3):237--252.

\bibitem[Yu, 1994]{yu94}
Yu, B. (1994).
\newblock Rates of convergence for empirical processes of stationary mixing
  sequences.
\newblock {\em The Annals of Probability}, Vol. 22,No. 1:94--116.

\bibitem[Zhou, 2012]{zhou2012}
Zhou, Z. (2012).
\newblock Measuring nonlinear dependence in time-series, a distance correlation
  approach.
\newblock {\em Journal of Time Series Analysis}, 33(3):438--457.

\end{thebibliography}

\newpage

\section{Appendix: technical details} \label{sec5}
\def\theequation{5.\arabic{equation}}
\setcounter{equation}{0}

The proofs of Theorems \ref{asympnormalityrho} and \ref{asympnormalitytau} rely on a blocking technique which will be summarized in Section \ref{blocking}. In Section \ref{sec:auc} we state and prove the covariance inequalities that are crucial in order to derive the convergence of the linear and degenerate parts of the \textit{U-lag-window estimate}. Finally, in Section \ref{sec:5.3} we provide the details for the proofs of results and equations given in Section \ref{sec4}.

\bigskip

For simplicity of notation, let $\bm{X}_{t,k}:=(X_{t},X_{t+k})^T$. 

\subsection{Blocking results for stationary $\beta$-mixing processes}\label{blocking}

In order to transfer classical results from the iid case to sums of $\beta$-mixing stationary time series, we apply a blocking technique with alternate "large" blocks of size $p_n$ and "small" blocks of size $q_n$ from \cite{arcones94} based on a blocking technique introduced by \cite{yu94} with blocks of equal size $p_n$. For each fixed $n$, we divide the original sequence $X:=(X_1,\dots,X_n)$ into $\mu_n$ blocks of size $p_n$ alternating with $\mu_n$ blocks of size $q_n$ and a remainder block $\mathcal{R}$ of length $n-2\mu_n$. The block size $q_n$ of the "small" $\Delta$ blocks is chosen depending on the mixing conditions on $X$ and the size of $r_n$. That is, $q_n$ is chosen large enough such that the $\Gamma$ blocks are "almost" independent of each other, but small enough such that the sequence composed of these $\Gamma$ blocks behaves similarly to the original mixing sequence. The block size $p_n$ is chosen analogously. More precisely, we assume that
\beq\label{mun}
(\mu_n-1)(p_n+q_n)<n\leq\mu_n(p_n+q_n)
\eeq
and define $\text{for }j=1,\dots,\mu_n$
\[
\ba{rcl}
\Gamma_j&:=&\{i:(j-1)(p_n+q_n)+1\leq i\leq (j-1)(p_n+q_n)+p_n\},\\
\Delta_j&:=&\{i:(j-1)(p_n+q_n)+p_n+1\leq i\leq j(p_n+q_n)\}\\
\Gamma_{\mathcal{R}}&:=&\{i:\mu_n(p_n+q_n)+1\leq i\leq n\wedge\mu_n(p_n+q_n)+p_n\}\\
\Delta_{\mathcal{R}}&:=&\{i:(n\wedge \Big[\mu_n(p_n+q_n)+p_n\Big])+1\leq i\leq n\}.\\
\mathcal{R}&:=&\Gamma_{\mathcal{R}}\cup\Delta_{\mathcal{R}}
\ea
\]
We denote the random variables of $X$ belonging to block $\Gamma_j$, $\Delta_j,\, j=1,\dots,\mu_n$ or $\mathcal{R}$ by
\[
X(\Gamma_j):=\{X_i:i\in\Gamma_j\},\quad X(\Delta_j):=\{X_i:i\in\Delta_j\},\quad X(\mathcal{R}):=\{X_i:i\in\mathcal{R}\},
\]
respectively. This yields a sequence of alternating $\Gamma$ and $\Delta$ blocks
\[
X(\Gamma_1),X(\Delta_1),X(\Gamma_2),X(\Delta_2),\dots,X(\Gamma_{\mu_n}),X(\Delta_{\mu_n}),X(\mathcal{R})
\]
We then construct a one-dependent sequence $Y$ of independent blocks defined as
\[
Y(\Gamma_1),Y(\Delta_1),Y(\Gamma_2),Y(\Delta_2),\dots,Y(\Gamma_{\mu_n}),Y(\Delta_{\mu_n}),Y(\mathcal{R}),
\]
and independent of the original sequence $X$. Furthermore, the blocks $Y(\Gamma_j):=\{Y_i:i\in\Gamma_j\}$, $Y(\Delta_j):=\{Y_i:i\in\Delta_j\},\, j=1,\dots,\mu_n$ and $Y(\mathcal{R}):=\{Y_i:i\in\mathcal{R}\}$ are identically distributed as the corresponding blocks in the sequence $X$, i.e.
\[
X(\Gamma_j)\overset{\mathcal{D}}{=}Y(\Gamma_j)\quad\text{and}\quad X(\Delta_j)\overset{\mathcal{D}}{=}Y(\Delta_j)\quad\text{and}\quad X(\mathcal{R})\overset{\mathcal{D}}{=}Y(\mathcal{R}) .
\]
The existence of a proper measurable space that hosts both sequences, $X$ and the independent block sequence $Y$, as well as measurability issues on this space are adressed in \cite{yu94}. Denote by $X_{\Gamma}$ and $Y_{\Gamma}$ block sequences corresponding to the $\Gamma$ blocks and by $X_{\Delta}$ and $Y_{\Delta}$ the block sequences corresponding to the $\Delta$ blocks, e.g.
\[
X_{\Gamma}:=X(\Gamma_1),X(\Gamma_2),\dots,X(\Gamma_{\mu_n})
\]
Note that we choose the block size $q_n$ such that the dependence between the blocks $X_{\Gamma}$ of the original $\beta$-mixing sequence $X$ becomes weaker as $q_n$ increases. The next lemma is a slightly adapted version of Lemma 4.1 in \cite{yu94} and is proven analogously. It shows that the $\Gamma$ or, respectively, $\Delta$ blocks of the original sequence $X$ can be related to the $\Gamma$ or, respectively, $\Delta$ blocks of the independent block sequence $Y$ in the following way.

\begin{lem}\label{link}
Denote by $Q$ and $\tilde Q$ be the distributions of $X_{\Gamma}$ and $Y_{\Gamma}$, respectively. Then, for any measurable function $g$ on $\R^{\mu_n p_n}$ with $\|g\|_{\infty}\leq M<\infty$,

\[
\Big|\E_{Q}[g(X_{\Gamma})]-\E_{\tilde Q}[g(Y_{\Gamma})]\Big|\leq M(\mu_n-1)\beta(q_n).
\]

Similarly, if $P$ and $\tilde P$ denote the distributions of $X_{\Delta}$ and $Y_{\Delta}$, respectively, and if $\tilde g$ is a measurable function on $\R^{\mu_n q_n}$ with $\|\tilde g\|_{\infty}\leq N<\infty$, then

\[
\Big|\E_{P}[\tilde g(X_{\Delta})]-\E_{\tilde P}[\tilde g(Y_{\Delta})]\Big|\leq N(\mu_n-1)\beta(p_n).
\]
\end{lem} 

In order to establish the convergence in probability of the parts of the \textit{U-lag-window estimate} corresponding to the linear and degenerate part in the Hoeffding decomposition we prove several covariance inequalities for $\beta$-mixing data. To this end we apply a coupling technique by \cite{berbee79}. The idea is to replace successively dependent variables by variables that have the same distribution but are independent of the original variables and all other involved variables with the smallest error possible. \cite{berbee79} found the following in the case of $\beta$-mixing data.

\begin{lem} [\cite{berbee79}]\label{berbee}
Suppose on a probability space there is defined a pair $(X,Y)$ of random variables with values in Borel spaces. If the probability space is rich enough, it can be extended with a random variable $Y^{\prime}$, independent of $X$ and distributed as $Y$ such that
\[
\Prob(Y^{\prime}\neq Y)=\frac{1}{2}\Big\|\Prob_{(X,Y)}-\Prob_X\otimes \Prob_Y\Big\|_{TV}=\beta(\sigma(X),\sigma(Y)).
\]
\end{lem}

\subsection{Auxiliary technical results} \label{sec:auc}

\begin{lem}\label{unifmom}
Let Assumption (C2) hold. Then, the kernels $h_{c,k}$ defined in \eqref{hkhoef}, $c=1,\dots,m$ have uniform $(2+\delta)$ moments, i.e. there exist $\delta, M_c>0$ such that for all $t_1,\dots,t_{c},k\in\Z$, $1\leq j\leq 2c$,
\al{
\max\Big\{\int_{\R}\dots\int_{\R}|h_{c,k}|^{2+\delta}dG_c,
\int_{\R}\dots\int_{\R}|h_{c,k}|^{2+\delta}dG_{j,c}^{(1)}dG_{j,c}^{(2)}\Big\}\leq M_0<\infty,
}
where $G_c$, $G_{j,c}^{(1)}$ and $G_{j,c}^{(2)}$ denote the joint distributions of $(X_{t_{(1)}},\dots,X_{t_{(2c)}})$, $(X_{t_{(1)}},\dots,X_{t_{(j)}})$ and  $(X_{t_{(j+1)}},\dots,X_{t_{(2c)}})$, respectively, with $(t_{(1)}, t_{(2)},\dots, t_{(2m)})$, $t_{(1)}\leq\dots\leq t_{(2c)}$ the sorted version of the vector $(t_1,t_1+k,t_2,t_2+k\,\dots,t_c,t_c+k)$.
\end{lem}

\begin{lem}\label{lemma1m}
Let $\bm{X}_{t_j,k}^*$ denote an independent and identically distributed copy of $\bm{X}_{t_j,k}$ on a possibly richer probability space that is independent of $\bm{X}_{t_1,k},\dots,\bm{X}_{t_{2c},k}$. Then, for $2\leq c\leq m$ and arbitrary $t_1,\dots,t_{2c}\in\Z$,
\[
\E\Big[h_{c,k}(\bm{X}_{t_1,k}^*,\bm{X}_{t_2,k},\dots,\bm{X}_{t_c,k})h_{c,k}(\bm{X}_{t_{c+1},k},\dots,\bm{X}_{t_{2c},k})\Big]=0,
\]
where the latter equation also holds if any other pair $\bm{X}_{t_j,k}$ is replaced by an iid copy $\bm{X}_{t_j,k}^*$.
\end{lem}

\begin{lem}\label{boundshm}
If Assumptions (C1) -- (C3) are satisfied, we have   for any fixed $0\leq k\leq \floor{r_n}$ 
\bi
\item[(1)] For any $t\in\Z$,
\als{\label{varbloc}
\Big|\E[h_{1,k}\Big(\bm{X}_{t,k}\Big)h_{1,k}\Big(\bm{X}_{t+l,k}\Big)]\Big|\leq \begin{cases}
2M_1^{\frac{2}{2+\delta}}\beta^{\frac{\delta}{2+\delta}}(l-k),&\text{ if }l>k,\\
8M_1^{\frac{2}{2+\delta}}\beta^{\frac{\delta}{2+\delta}}(\min\{l,k-l\}),&\text{ if }0\leq l\leq k.
\end{cases}
}

\item[(2)] If  $1\leq t_1<t_2<\dots <t_{2c}\leq n-k$ and
\al{
m(t_1,\dots,t_{2c})
:={}&\max\big\{\min\{t_2-t_1,(t_1+k)-t_2\},\min\{t_{2c}-t_{2c-1},(t_{2c-1}+k)-t_{2c}\}\big\}
}
we have for any permutation
 $\gamma$  of $\{1,\dots,2c\}$
\bi
\item[(i)]
~~$
|\E[h_{c,k}(\bm{X}_{t_{\gamma(1)},k},\dots,\bm{X}_{t_{\gamma(c)},k})h_{c,k}(\bm{X}_{t_{\gamma(c+1)},k},\dots,\bm{X}_{t_{\gamma(2c)},k})]|\leq M_c^{\frac{2}{2+\delta}}.
$
\item[(ii)] if $(t_2-t_1)>k$ or $(t_{2c}-t_{2c-1})>k$,
\als{
&|\E[h_{c,k}(\bm{X}_{t_{\gamma(1)},k},\dots,\bm{X}_{t_{\gamma(c)},k})h_{c,k}(\bm{X}_{t_{\gamma(c+1)},k},\dots,\bm{X}_{t_{\gamma(2c)},k})]|\notag\\
\leq & 2M_c^{\frac{2}{2+\delta}}\beta^{\frac{\delta}{2+\delta}}(\max\Big\{t_2-(t_1+k),t_{2c}-(t_{2c-1}+k)\Big\}).\label{faraway3}
}
\item[(iii)] if $t_2-t_1\leq k$, $t_{2c}-t_{2c-1}\leq k$, $t_3-t_2>2k$ and $t_{2c-1}-t_{2(c-1)}>2k$, then
\als{
&|\E[h_{c,k}(\bm{X}_{t_{\gamma(1)},k},\dots,\bm{X}_{t_{\gamma(c)},k})h_{c,k}(\bm{X}_{t_{\gamma(c+1)},k},\dots,\bm{X}_{t_{\gamma(2c)},k})]|
\leq{} 12M_c^{\frac{2}{2+\delta}}\beta^{\frac{\delta}{2+\delta}}(m).\label{close3}
}

\ei
\ei
\end{lem}

\subsubsection{Proof of Lemma~\ref{unifmom} }
From assumption (C2), we have
\al{
\E\Big|\E[h(\bm{Y}^{(1)},\dots,\bm{Y}^{(m)})|\bm{Y}^{(1)},\dots,\bm{Y}^{(c)}]-\E[h(\bm{Y}^{(1)},\dots,\bm{Y}^{(m)})]\Big|^{2+\delta}
\leq{}&2^{2+\delta}M_0
}
and therefore, by the definition of the Hoeffding decomposition,
\[
\E\Big|h_{1,k}\Big(\bm{Y_1}\Big)\Big|^{2+\delta}\leq{}2^{2+\delta}M_0.
\]
As $h_{c,k}$ is recursively defined by
\al{
h_{c,k}(\bm{y_1},\dots,\bm{y_c})
={}&\E[h(\bm{Y}^{(1)},\dots,\bm{Y}^{(m)})|\bm{Y}^{(1)}=\bm{y_1},\dots,\bm{Y}^{(c)}=\bm{y_c}]-\E[h(\bm{Y}^{(1)},\dots,\bm{Y}^{(m)})]\\
&-\sum_{j=1}^{c-1}\sum\limits_{\substack{\{\nu_1,\dots,\nu_j\}\subset\{1,\dots,c\}\\ \nu_1<\dots<\nu_j}}h_{j,k}(\bm{y_{\nu_1}},\dots,\bm{y_{\nu_j}}),
}
$h_{c,k}$ also has uniform $(2+\delta)$-moments. \hfill $\Box$
 

\subsubsection{Proof of Lemma \ref{lemma1m}} Recall the following property of the conditional expectation [see Theorem 6.4 in \cite{kallenberg10}] which can easily be adapted to more than one $\mathcal{F}$-measurable random variable:\\

\textit{Let $X$ and $Y$ be random variables and $\mathcal{F}$ a $\sigma$-algebra such that $X$ is $\mathcal{F}$-measurable and $Y$ is independent of $\mathcal{F}$. Then, for any measurable function $f(x,y)$ with $\E|f(X,Y)|<\infty$},

\be\label{condexp}
\E[f(X,Y)|\mathcal{F}]=F(X) \qquad \text{a.s.},
\ee

\textit{where $F(x)=\E[f(x,Y)]$.}\\

Thus, by the law of total expectation, we have
\al{
&\E\Big[h_{c,k}(\bm{X}^*_{t_1,k},\bm{X}_{t_2,k},\dots,\bm{X}_{t_c,k})h_{c,k}(\bm{X}_{t_{c+1},k},\dots,\bm{X}_{t_{2c},k})\Big]\\
={}&\E\Big[\E\Big[h_{c,k}(\bm{X}^*_{t_1,k},\bm{X}_{t_2,k},\dots,\bm{X}_{t_c,k})h_{c,k}(\bm{X}_{t_{c+1},k},\dots,\bm{X}_{t_{2c},k})|\bm{X}_{t_2,k},\dots,\bm{X}_{t_{2c},k}\Big]\Big]\\
={}&\E\Big[\E\Big[h_{c,k}(\bm{X}^*_{t_1,k},\bm{X}_{t_2,k},\dots,\bm{X}_{t_c,k})|\bm{X}_{t_2,k},\dots,\bm{X}_{t_{2c},k}\Big]h_{c,k}(\bm{X}_{t_{c+1},k},\dots,\bm{X}_{t_{2c},k})\Big]
}

Obviously, $\bm{X}_{t_2,k},\dots,\bm{X}_{t_c,k}$ are $\sigma(\bm{X}_{t_2,k},\dots,\bm{X}_{t_{2c},k})$-measurable and $\bm{X}^*_{t_1,k}$ is independent of \\ $\sigma(\bm{X}_{t_2,k},\dots,\bm{X}_{t_{2c},k})$. As, additionally, $\E|h_{c,k}(\bm{X}^*_{t_1,k},\bm{X}_{t_2,k},\dots,\bm{X}_{t_c,k})|<\infty$ by Lemma \ref{unifmom} we have that 
\al{
\E[h_{c,k}(&\bm{X}^*_{t_1,k},\bm{X}_{t_2,k},\dots,\bm{X}_{t_c,k})|\bm{X}_{t_2,k},\dots,\bm{X}_{t_{2c},k}]=:H_{c-1,k}(\bm{X}_{t_2,k},\dots,\bm{X}_{t_{c},k}),
}
where
\[
H_{c-1,k}(\bm{y_2},\dots,\bm{y_c})=\E[h_{c,k}(\bm{X}^*_{t_1,k},\bm{y_2},\dots,\bm{y_c})].
\]
We will now show that $H_{c-1,k}(\bm{y_2},\dots,\bm{y_c})=0$. To this end, we consider the integral representation of $h_{c,k}$, i.e. similarly as in the proof of Theorem 2 in [\cite{lee90}, pg. 28], we obtain by the symmetry of $h$,
\al{
h_{c,k}(\bm{y_1},\dots,\bm{y_c})
={}&\int_{\R^2}\cdots\int_{\R^2}h(\bm{u_1},\dots,\bm{u_m})\prod_{j=1}^c(dG_{\bm{y_j}}(\bm{u_j})-dF_k(\bm{u_j}))\prod_{j=c+1}^{m}dF_k(\bm{u_j})\\
={}&\int_{\R^2}\cdots\int_{\R^2}h(\bm{y_1},\bm{u_2},\dots,\bm{u_m})\prod_{j=2}^c(dG_{\bm{y_j}}(\bm{u_j})-dF_k(\bm{u_j}))\prod_{j=c+1}^{m}dF_k(\bm{u_j})\\
&-\int_{\R^2}\cdots\int_{\R^2}h(\bm{u_1},\dots,\bm{u_m})dF_k(\bm{u_1})\prod_{j=2}^c(dG_{\bm{y_j}}(\bm{u_j})-dF_k(\bm{u_j}))\prod_{j=c+1}^{m}dF_k(\bm{u_j})\\
={}&\int_{\R^2}\cdots\int_{\R^2}h(\bm{y_1},\bm{u_2}\dots,\bm{u_m})\prod_{j=2}^c(dG_{\bm{y_j}}(\bm{u_j})-dF_k(\bm{u_j}))\prod_{j=c+1}^{m}dF_k(\bm{u_j})\\
&-h_{c-1,k}(\bm{y_2},\dots,\bm{y_c}).
}
Integrating both sides with respect to $\bm{u_1}\sim F_k$ yields
\al{
&\int_{\R^{2}}h_{c,k}(\bm{u_1},\bm{y_2},\dots,\bm{y_c})dF_k(\bm{u_1})\\
={}&\int_{\R^2}\cdots\int_{\R^2}h(\bm{u_1},\dots,\bm{u_m})dF_k(\bm{u_1})\prod_{j=2}^c(dG_{\bm{y_j}}(\bm{u_j})-dF_k(\bm{u_j}))\prod_{j=c+1}^{m}dF_k(\bm{u_j})\\
&-h_{c-1,k}(\bm{y_2},\dots,\bm{y_c})\\
={}&h_{c-1,k}(\bm{y_2},\dots,\bm{y_c})-h_{c-1,k}(\bm{y_2},\dots,\bm{y_c})={}0.
}
Observing that 
\al{
\int_{\R^{2}}h_{c,k}(\bm{u_1},\bm{y_2},\dots,\bm{y_c})dF_k(\bm{u_1})=\E[h_{c,k}(\bm{X}^*_{t_1,k},\bm{y_2},\dots,\bm{y_c})]=H_{c-1,k}(\bm{y_2},\dots,\bm{y_c})
}
we have
\al{
\E\Big[h_{c,k}(\bm{X}^*_{t_1,k},\bm{X}_{t_2,k},\dots,\bm{X}_{t_c,k})|\bm{X}_{t_2,k},\dots,\bm{X}_{t_{2c},k}\Big]=H_{c-1,k}(\bm{X}_{t_2,k},\dots,\bm{X}_{t_{c},k})=0
}
and altogether,
\al{
&\E\Big[h_{c,k}(\bm{X}^*_{t_1,k},\bm{X}_{t_2,k},\dots,\bm{X}_{t_c,k})h_{c,k}(\bm{X}_{t_{c+1},k},\dots,\bm{X}_{t_{2c},k})\Big]\\
={}&\E\Big[\E\Big[h_{c,k}(\bm{X}^*_{t_1,k},\bm{X}_{t_2,k},\dots,\bm{X}_{t_c,k})|\bm{X}_{t_2,k},\dots,\bm{X}_{t_{2c},k}\Big]h_{c,k}(\bm{X}_{t_{c+1},k},\dots,\bm{X}_{t_{2c},k})\Big]={}0
}
which concludes the proof. \hfill $\Box$


\subsubsection{Proof of Lemma \ref{boundshm}}\label{sec:proofboundshm}

\bi[leftmargin=*]
\item[(1)] If $l>k\geq 0$ replace the pair $\bm{X}_{t,k}:=\colvec{X_t}{X_{t+k}}$ using Berbee's coupling technique by an identically distributed copy $\bm{X}_{t,k}^*
$ that is independent of $\bm{X}_{t,k}$ and $\bm{X}_{t+l,k}$ and such that 
\[
\Prob(\bm{X}_{t,k}\neq \bm{X}_{t,k}^*)\leq \frac{1}{2}\Big\|\Prob_{(\bm{X}_{t,k},\bm{X}_{t+l,k})}-\Prob_{\bm{X}_{t,k}}\otimes\Prob_{\bm{X}_{t+l,k}}\Big\|_{TV}=\beta(l-k).
\]
Then, 
\bea
&&\Big|\E[h_{1,k}(\bm{X}_{t,k})h_{1,k}(\bm{X}_{t+l,k})]-\E[h_{1,k}(\bm{X}_{t,k}^*)h_{1,k}(\bm{X}_{t+l,k})]\Big|\\
&\leq &\Big|\E\Big[\Big(h_{1,k}(\bm{X}_{t,k})-h_{1,k}(\bm{X}_{t,k}^*)\Big)h_{1,k}(\bm{X}_{t+l,k})I(\bm{X}_{t,k}\neq \bm{X}_{t,k}^*)\Big]\Big|\\
\eea
By H\"older's inequality $(\frac{1}{2+\delta}+\frac{1}{2+\delta}+\frac{1}{\frac{\delta+2}{\delta}}=1)$ we obtain
\begin{multline*}
\Big|\E[\Big(h_{1,k}(\bm{X}_{t,k})-h_{1,k}(\bm{X}_{t,k}^*)\Big)h_{1,k}(\bm{X}_{t+l,k})I(\bm{X}_{t,k}\neq \bm{X}_{t,k}^*)]\Big|\\
\leq{} 2M_1^{\frac{2}{2+\delta}}\Big(\E[I(\bm{X}_{t,k}\neq \bm{X}_{t,k}^*)]\Big)^{\frac{\delta}{2+\delta}}
\leq{} 2M_1^{\frac{2}{2+\delta}}\beta^{\frac{\delta}{2+\delta}}(l-k),
\end{multline*}
which gives
\al{
\Big|\E[h_{1,k}(\bm{X}_{t,k})h_{1,k}(\bm{X}_{t+l,k})]\Big|\leq 2M_1^{\frac{2}{2+\delta}}\beta^{\frac{\delta}{2+\delta}}(l-k)+\Big|\E[h_{1,k}(\bm{X}_{t,k}^*)h_{1,k}(\bm{X}_{t+l,k})]\Big|.
}
Now, $h_{1,k}(\bm{X_{t}^*})$ is independent of $h_{1,k}(\bm{X_{t+l}})$ and the result follows by Lemma \ref{lemma1m}.\\

If $0\leq l\leq k$, using Berbee's coupling technique, first replace $X_t$ by an identically distributed copy $X_{t}^{*}$ that is independent of $X_{t+k}$, $X_{t+l}$ and $X_{t+l+k}$ and such that
\[
\Prob(X_{t}^{*}\neq X_{t})\leq \frac{1}{2}\Big\|\Prob_{(X_t,X_{t+l},X_{t+k},X_{t+l+k})}-\Prob_{X_t}\otimes\Prob_{(X_{t+l},X_{t+k},X_{t+l+k})}\Big\|_{TV}=\beta(l).
\]
By H\"older's inequality, we obtain with similar arguments as in the case $l>k\geq 0$ that
\[
\Big|\E\Big[\hun{X_{t}}{X_{t+k}}\hun{X_{t+l}}{X_{t+k+l}}\Big]\Big|\leq 2M_1^{\frac{2}{2+\delta}}\beta^{\frac{\delta}{2+\delta}}(l)+\Big|\E\Big[\hun{X_{t}^{*}}{X_{t+k}}\hun{X_{t+l}}{X_{t+k+l}}\Big]\Big|,
\]
Then, replace $X_{t+l}$ by an independent copy $X_{t+l}^*$ that is independent of $X_{t}^{*}$, $X_{t+k}$, $X_{t+l}$ and $X_{t+l+k}$ such that
\al{
\Prob(X_{t+l}^{*}\neq X_{t+l})\leq &\frac{1}{2}\Big\|\Prob_{(X_t^*,X_{t+l},X_{t+k},X_{t+l+k})}-\Prob_{X_{t+l}}\otimes\Prob_{(X_t^*,X_{t+k},X_{t+l+k})}\Big\|_{TV}\\
\leq{}&\frac{1}{2}\Big\|\Prob_{X_t^*}\otimes\Prob_{(X_{t+l},X_{t+k},X_{t+l+k})}-\Prob_{X_t^*}\otimes\Prob_{X_{t+l}}\otimes\Prob_{(X_{t+k},X_{t+l+k})}\Big\|_{TV}\\
\leq{}&\frac{1}{2}\Big\|\Prob_{(X_{t+l},X_{t+k},X_{t+l+k})}-\Prob_{X_{t+l}}\otimes\Prob_{(X_{t+k},X_{t+l+k})}\Big\|_{TV}\\
={}&\beta(k-l),
}
where we have used Lemma 1 in \cite{eberlein84}. Then,
\[
\Big|\E\Big[\hun{X_{t}^*}{X_{t+k}}\hun{X_{t+l}}{X_{t+k+l}}\Big]\Big|\leq 2M_1^{\frac{2}{2+\delta}}\beta^{\frac{\delta}{2+\delta}}(k-l)+\Big|\E\Big[\hun{X_{t}^{*}}{X_{t+k}}\hun{X_{t+l}^*}{X_{t+k+l}}\Big]\Big|.
\]
Finally, replace $X_{t+k}$ by an independent copy $X_{t+k}^*$ that is independent of $X_{t}^{*}$, $X_{t+l}^*$, $X_{t+k}$ and $X_{t+l+k}$ and such that 
\al{
\Prob(X_{t+k}^{*}\neq X_{t+k})\leq &\frac{1}{2}\Big\|\Prob_{(X_t^*,X_{t+l}^*,X_{t+k},X_{t+l+k})}-\Prob_{X_{t+k}}\otimes\Prob_{(X_t^*,X_{t+l}^*,X_{t+l+k})}\Big\|_{TV}\\
\leq{}&\frac{1}{2}\Big\|\Prob_{X_t^*}\otimes\Prob_{X_{t+l}^*}\otimes\Prob_{(X_{t+k},X_{t+l+k})}-\Prob_{X_t^*}\otimes\Prob_{X_{t+l}^*}\otimes\Prob_{X_{t+k}}\otimes\Prob_{X_{t+l+k}}\Big\|_{TV}\\
\leq{}&\frac{1}{2}\Big\|\Prob_{(X_{t+k},X_{t+l+k})}-\Prob_{X_{t+k}}\otimes\Prob_{X_{t+l+k}}\Big\|_{TV}\\
={}&\beta(l)
}
which gives
\[
\Big|\E\Big[\hun{X_{t}^*}{X_{t+k}}\hun{X_{t+l}^*}{X_{t+k+l}}\Big]\Big|\leq 2M_1^{\frac{2}{2+\delta}}\beta^{\frac{\delta}{2+\delta}}(l)+\Big|\E\Big[\hun{X_{t}^{*}}{X_{t+k}^*}\hun{X_{t+l}^*}{X_{t+k+l}}\Big]\Big|.
\]
Altogether,
\al{
\Big|\E\Big[\hun{X_{t}}{X_{t+k}}\hun{X_{t+l}}{X_{t+k+l}}\Big]\Big|\leq{}& 6M_1^{\frac{2}{2+\delta}}\beta^{\frac{\delta}{2+\delta}}(\min\Big\{l,k-l\Big\})\\
&+\Big|\E\Big[\hun{X_{t}^{*}}{X_{t+k}^*}\Big]\E\Big[\hun{X_{t+l}^*}{X_{t+k+l}}\Big]\Big|.
}
Next, observe that the last summand does not vanish as $\colvec{X_{t}^{*}}{X_{t+k}^*}$ does not have distribution $F_k$. Therefore, using Berbee's coupling technique, we rereplace $X_{t}^{*}$ by an independent copy $X_{t}^{\circ}$ such that the couple $\colvec{X_{t}^{\circ}}{X_{t+k}^*}$ has distribution $F_k$, is independent of $\colvec{X_{t+l}^*}{X_{t+k+l}}$ and
\al{
\Prob(X_{t}^{\circ}\neq X_{t}^*)\leq{}&\frac{1}{2}\Big\|\Prob_{X_t^*}\otimes\Prob_{X_{t+k}^*}-\Prob_{(X_t^{\circ},X_{t+k}^*)}\Big\|_{TV}=\beta(k).
}
Hence,
\al{
\Big|\E\Big[\hun{X_{t}^{*}}{X_{t+k}^*}-\hun{X_{t}^{\circ}}{X_{t+k}^*}\Big]\Big|\leq 2M_1^{\frac{1}{2+\delta}}\beta^{\frac{\delta}{2+\delta}}(k). 
}
Observing that $\min\Big\{l,k-l\Big\}\leq k$ and that the $\beta$-mixing coefficients are monotone decreasing, we have
\[
\beta^{\frac{\delta}{2+\delta}}(k)\leq \beta^{\frac{\delta}{2+\delta}}(\min\Big\{l,k-l\Big\})
\]
and as $\E\Big[\hun{X_{t}^{\circ}}{X_{t+k}^*}\Big]=0$, we can conclude that
\[
\Big|\E\Big[\hun{X_{t}}{X_{t+k}}\hun{X_{t+l}}{X_{t+k+l}}\Big]\Big|\leq 8M_1^{\frac{2}{2+\delta}}\beta^{\frac{\delta}{2+\delta}}(\min\Big\{l,k-l\Big\}).
\]
\item[(2) ]
\bi[leftmargin=*]
\item[\textbf{(i)}] By H\"older's inequality $(\frac{1}{2+\delta}+\frac{1}{\frac{2+\delta}{1+\delta}}=1)$ and as $\frac{2+\delta}{1+\delta}< 2+\delta$ we obtain
\bea
&&|\E[h_{c,k}(\bm{X}_{t_{\gamma(1)},k},\dots,\bm{X}_{t_{\gamma(c)},k})h_{c,k}(\bm{X}_{t_{\gamma(c+1)},k},\dots,\bm{X}_{t_{\gamma(2c)},k})]|\\
&\leq &\Big(\E|h_{c,k}(\bm{X}_{t_{\gamma(1)},k},\dots,\bm{X}_{t_{\gamma(c)},k})|^{2+\delta}\Big)^{\frac{1}{2+\delta}}\Big(\E|h_{c,k}(\bm{X}_{t_{\gamma(c+1)},k},\dots,\bm{X}_{t_{\gamma(2c)},k})|^{\frac{2+\delta}{1+\delta}}\Big)^{\frac{1+\delta}{2+\delta}}\\
&\leq &\Big(\E|h_{c,k}(\bm{X}_{t_{\gamma(1)},k},\dots,\bm{X}_{t_{\gamma(c)},k})|^{2+\delta}\Big)^{\frac{1}{2+\delta}}\Big(\E|h_{c,k}(\bm{X}_{t_{\gamma(c+1)},k},\dots,\bm{X}_{t_{\gamma(2c)},k})|^{2+\delta}\Big)^{\frac{1}{2+\delta}}.
\eea
where we have used that $(\E[|Z|^p])^{1/p}\leq (\E[|Z|^q])^{1/q}$ for $0<p\leq q$. Hence, by Lemma \ref{unifmom},
\bea
|\E[h_{c,k}(\bm{X}_{t_{\gamma(1)},k},\dots,\bm{X}_{t_{\gamma(c)},k})h_{c,k}(\bm{X}_{t_{\gamma(c+1)},k},\dots,\bm{X}_{t_{\gamma(2c)},k})]|&\leq &M_c^{\frac{2}{2+\delta}}.
\eea
\item[\textbf{(ii)}] For brevity, we only consider the case $\gamma=\textrm{id}$. The other cases are treated similarly but require a more complex notation.\\
In order to prove inequality (\ref{faraway3}), according to the coupling Lemma \ref{berbee} by \cite{berbee79}, depending on whether $(t_2-t_1)>(t_{2c}-t_{2c-1})>k$ or $k<(t_2-t_1)\leq(t_{2c}-t_{2c-1})$, we can choose  a random variable $\bm{X}^*_{t_1,k}$ or respectively $\bm{X}^*_{t_{2c},k}$ that has the same distribution as $\bm{X}_{t_1,k}$ or respectively $\bm{X}_{t_{2c},k}$, independent of $\bm{X}_{t_1,k},\dots,\bm{X}_{t_{2c},k}$ and such that
\[
\Prob(\bm{X}^*_{t_{1},k}\neq\bm{X}_{t_{1},k})\leq \beta(t_2-(t_1+k))
\] 
First, consider the case where $(t_2-t_1)>(t_{2c}-t_{2c-1})>k$, that is we replace $\bm{X}_{t_1,k}$ by an independent identically distributed copy $\bm{X}^*_{t_1,k}$. Then, by Lemma \ref{lemma1m},
\al{
&|\E[h_{c,k}(\bm{X}_{t_{1},k},\dots,\bm{X}_{t_{c},k})h_{c,k}(\bm{X}_{t_{c+1},k},\dots,\bm{X}_{t_{2c},k})]|\\
={}&|\E[(h_{c,k}(\bm{X}_{t_{1},k},\dots,\bm{X}_{t_{c},k})-h_{c,k}(\bm{X}^*_{t_{1},k},\bm{X}_{t_{2},k},\dots,\bm{X}_{t_{c},k}))h_{c,k}(\bm{X}_{t_{c+1},k},\dots,\bm{X}_{t_{2c},k})]|
}
Splitting the probability space, we obtain
\al{
&|\E[(h_{c,k}(\bm{X}_{t_{1},k},\dots,\bm{X}_{t_{c},k})-h_{c,k}(\bm{X}^*_{t_{1},k},\bm{X}_{t_{2},k},\dots,\bm{X}_{t_{c},k}))h_{c,k}(\bm{X}_{t_{c+1},k},\dots,\bm{X}_{t_{2c},k})]|\\
&\leq\E[|(h_{c,k}(\bm{X}_{t_{1},k},\dots,\bm{X}_{t_{c},k})-h_{c,k}(\bm{X}^*_{t_{1},k},\bm{X}_{t_{2},k},\dots,\bm{X}_{t_{c},k}))\\
&\hspace{0.5cm}\cdot h_{c,k}(\bm{X}_{t_{c+1},k},\dots,\bm{X}_{t_{2c},k})I(\bm{X}^*_{t_{1},k}\neq\bm{X}_{t_{1},k})|]\\
&+\E[|(h_{c,k}(\bm{X}_{t_{1},k},\dots,\bm{X}_{t_{c},k})-h_{c,k}(\bm{X}^*_{t_{1},k},\bm{X}_{t_{2},k},\dots,\bm{X}_{t_{c},k}))\\
&\hspace{0.5cm}\cdot h_{c,k}(\bm{X}_{t_{c+1},k},\dots,\bm{X}_{t_{2c},k})I(\bm{X}^*_{t_{1},k}=\bm{X}_{t_{1},k})|].
}

The second summand vanishes and for the first summand H\"older's inequality $(\frac{1}{2+\delta}+\frac{1}{2+\delta}+\frac{1}{\frac{2+\delta}{\delta}}=1)$ yields

\als{
&\E[|(h_{c,k}(\bm{X}_{t_{1},k},\dots,\bm{X}_{t_{c},k})-h_{c,k}(\bm{X}^*_{t_{1},k},\bm{X}_{t_{2},k},\dots,\bm{X}_{t_{c},k}))\notag\\
&\hspace{0.5cm}\cdot h_{c,k}(\bm{X}_{t_{c+1},k},\dots,\bm{X}_{t_{2c},k})I(\bm{X}^*_{t_{1},k}\neq\bm{X}_{t_{1},k})|]\notag\\
\leq{} &(\E|h_{c,k}(\bm{X}_{t_{1},k},\dots,\bm{X}_{t_{c},k})-h_{c,k}(\bm{X}^*_{t_{1},k},\bm{X}_{t_{2},k},\dots,\bm{X}_{t_{c},k})|^{2+\delta})^{\frac{1}{2+\delta}}\notag\\
&\cdot(\E|h_{c,k}(\bm{X}_{t_{c+1},k},\dots,\bm{X}_{t_{2c},k})|^{2+\delta})^{\frac{1}{2+\delta}}(\Prob(\bm{X}^*_{t_{1},k}\neq\bm{X}_{t_{1},k}))^{\frac{\delta}{2+\delta}}\notag\\
\leq{}&2M_c^{\frac{2}{2+\delta}}\beta^{\frac{\delta}{2+\delta}}(t_2-(t_1+k)),\label{unm}
}

where the latter inequality is due to Lemma \ref{berbee}. In the case where $k<(t_2-t_1)\leq(t_{2c}-t_{2c-1})$, we obtain

\als{\label{deuxm}
&|\E[(h_{c,k}(\bm{X}_{t_{1},k},\dots,\bm{X}_{t_{c},k})-h_{c,k}(\bm{X}^*_{t_{1},k},\bm{X}_{t_{2},k},\dots,\bm{X}_{t_{c},k}))h_{c,k}(\bm{X}_{t_{c+1},k},\dots,\bm{X}_{t_{2c},k})]|\notag\\
&\leq{} 2M_c^{\frac{2}{2+\delta}}\beta^{\frac{\delta}{2+\delta}}(t_{2c}-(t_{2c-1}+k))
}

Inequalities (\ref{unm}) and (\ref{deuxm}) together yield result (ii).

\item[\textbf{(iii)}] As in (ii), we only consider the case $\gamma=\textrm{id}$. Then, if $$\min\{t_2-t_1,(t_1+k)-t_2\}\geq \min\{t_{2c}-t_{2c-1},(t_{c-1}+k)-t_{2c}\}$$ replace one after another $X_{t_1}$, $X_{t_2}$, $X_{t_1+k}$ and $X_{t_2+k}$ according to Lemma \ref{berbee} by independent identically distributed copies $X_{t_1}^{\prime}$, $X_{t_2}^{\prime}$, $X_{t_1+k}^{\prime}$ and $X_{t_2+k}^{\prime}$ that are independent of the other involved random variables. Denote by $\bm{X}^{\prime}_{t_{j},k}$ the pair $(X_{t_j}^{\prime},X_{t_j+k}^{\prime})^T$, where $j=1,2$. Then, by Lemma \ref{unifmom} and similarly as in the proof of part (i)

\als{\label{betaminm}
&|\E[(h_{c,k}(\bm{X}_{t_{1},k},\dots,\bm{X}_{t_{c},k})-h_{c,k}(\bm{X}^{\prime}_{t_{1},k},\bm{X}^{\prime}_{t_{2},k},\bm{X}_{t_{3},k},\dots,\bm{X}_{t_{c},k}))h_{c,k}(\bm{X}_{t_{c+1},k},\dots,\bm{X}_{t_{2c},k})]|\notag\\
&\leq 8M_c^{\frac{2}{2+\delta}}\beta^{\frac{\delta}{2+\delta}}(\min\{t_2-t_1,(t_1+k)-t_2,t_3-(t_2+k)\})\notag\\
&=8M_c^{\frac{2}{2+\delta}}\beta^{\frac{\delta}{2+\delta}}(\min\{t_2-t_1,(t_1+k)-t_2\}),
}

where the latter equality is due to the assumption that $(t_3-t_2)>2k$. Hence,

\al{
|\E[h_{c,k}(\bm{X}_{t_{1},k}&,\dots,\bm{X}_{t_{c},k})h_{c,k}(\bm{X}_{t_{c+1},k},\dots,\bm{X}_{t_{2c},k})]|\\
\leq{}& 8M_c^{\frac{2}{2+\delta}}\beta^{\frac{\delta}{2+\delta}}(\min\{t_2-t_1,(t_1+k)-t_2\})\\
&+|\E[h_{c,k}(\bm{X}^{\prime}_{t_{1},k},\bm{X}^{\prime}_{t_{2},k},\bm{X}_{t_{3},k},\dots,\bm{X}_{t_{c},k})h_{c,k}(\bm{X}_{t_{c+1},k},\dots,\bm{X}_{t_{2c},k})]|
}

Note that the second summand does not necessarily vanish because, having replaced $X_{t_1}$, $X_{t_1+k}$, $X_{t_2}$ and $X_{t_2+k}$ by independent identically distributed copies one after another, the couple $(X_{t_j}^{\prime},X_{t_j+k}^{\prime})^T$ does not have distribution $F_k$. However, it is possible to bound $$|\E[h_{c,k}(\bm{X}^{\prime}_{t_{1},k},\bm{X}^{\prime}_{t_{2},k},\bm{X}_{t_{3},k},\dots,\bm{X}_{t_{c},k})h_{c,k}(\bm{X}_{t_{c+1},k},\dots,\bm{X}_{t_{2c},k})]|$$ by applying Lemma \ref{berbee} again in order to "{}rereplace"{} successively $\bm{X}^{\prime}_{t_{1},k}$ and $\bm{X}^{\prime}_{t_{2},k}$ by independent pairs $\bm{X}^{\circ}_{t_{1},k}$ and $\bm{X}^{\circ}_{t_{2},k}$ with distribution $F_k$. More precisely, according to Lemma \ref{berbee} we can replace $X_{t_1}^{\prime}$ by a random variable $X_{t_1}^{\circ}$ with the same distribution as $X_{t_1}^{\prime}$ that is independent of the other involved variables such that the couple $\bm{X}^{\circ}_{t_{1},k}:=(X_{t_1}^{\circ},X_{t_1+k}^{\prime})^T$ has distribution $F_k$. Then, similarly as in the proof of (ii),

\al{
&|\E[h_{c,k}(\bm{X}^{\prime}_{t_{1},k},\bm{X}^{\prime}_{t_{2},k},\bm{X}_{t_{3},k},\dots,\bm{X}_{t_{c},k})h_{c,k}(\bm{X}_{t_{c+1},k},\dots,\bm{X}_{t_{2c},k})]|\\
\leq{}& 2M_c^{\frac{1}{2+\delta}}\beta^{\frac{\delta}{2+\delta}}(k)+|\E[h_{c,k}(\bm{X}^{\circ}_{t_{1},k},\bm{X}^{\prime}_{t_{2},k},\bm{X}_{t_{3},k},\dots,\bm{X}_{t_{c},k})h_{c,k}(\bm{X}_{t_{c+1},k},\dots,\bm{X}_{t_{2c},k})]|,
}

where we have used that due to the independence of $X_{t_1}^{\prime}$, $X_{t_1+k}^{\prime}$ and all other involved variables, $\Prob(X_{t_1}^{\prime}\neq X_{t_1}^{\circ})=\beta(k)$. Analogously, we can replace $X_{t_2}^{\prime}$ by a random variable $X_{t_2}^{\circ}$ such that the couple $\bm{X}^{\circ}_{t_{2},k}:=(X_{t_2}^{\circ},X_{t_2+k}^{\prime})^T$ has distribution $F_k$ and is independent of $\bm{X}^{\circ}_{t_{1},k}$. Then,

\al{
&|\E[h_{c,k}(\bm{X}^{\prime}_{t_{1},k},\bm{X}^{\prime}_{t_{2},k},\bm{X}_{t_{3},k},\dots,\bm{X}_{t_{c},k})h_{c,k}(\bm{X}_{t_{c+1},k},\dots,\bm{X}_{t_{2c},k})]|\\
\leq{}& 2M_c^{\frac{1}{2+\delta}}\beta^{\frac{\delta}{2+\delta}}(k)+|\E[h_{c,k}(\bm{X}^{\circ}_{t_{1},k},\bm{X}^{\circ}_{t_{2},k},\bm{X}_{t_{3},k},\dots,\bm{X}_{t_{c},k})h_{c,k}(\bm{X}_{t_{c+1},k},\dots,\bm{X}_{t_{2c},k})]|,
}

Then, $\bm{X}^{\circ}_{t_{1},k}$ and $\bm{X}^{\circ}_{t_{2},k}$ both have distribution $F_k$ and are independent. Similar arguments as in the proof of Lemma \ref{lemma1m} (ii) show that also $$\E[h_{c,k}(\bm{X}^{\circ}_{t_{1},k},\bm{X}^{\circ}_{t_{2},k},\bm{X}_{t_{3},k},\dots,\bm{X}_{t_{c},k})h_{c,k}(\bm{X}_{t_{c+1},k},\dots,\bm{X}_{t_{2c},k})]=0$$ and hence,

\beq\label{boundjm}
|\E[h_{c,k}(\bm{X}^{\prime}_{t_{1},k},\bm{X}^{\prime}_{t_{2},k},\bm{X}_{t_{3},k},\dots,\bm{X}_{t_{c},k})h_{c,k}(\bm{X}_{t_{c+1},k},\dots,\bm{X}_{t_{2c},k})]|\leq 4M_c^{\frac{1}{2+\delta}}\beta^{\frac{\delta}{2+\delta}}(k).
\eeq

Observing that by the definition of the $\beta$-mixing coefficients and since $(t_2-t_1)\leq k$

\[
\beta^{\frac{\delta}{2+\delta}}(k)\leq \beta^{\frac{\delta}{2+\delta}}(\min\{t_2-t_1,(t_1+k)-t_2\}),
\]

equations (\ref{betaminm}) and (\ref{boundjm}) yield 

\al{
|\E[h_{c,k}(\bm{X}_{t_{1},k},\dots,\bm{X}_{t_{c},k})h_{c,k}(\bm{X}_{t_{c+1},k},\dots,\bm{X}_{t_{2c},k})]|
\leq{}& 12M_c^{\frac{2}{2+\delta}}\beta^{\frac{\delta}{2+\delta}}(\min\{t_2-t_1,(t_1+k)-t_2\}).
}

Analogously, if $\min\{t_2-t_1,(t_1+k)-t_2\} < \min\{t_{2c}-t_{2c-1},(t_{2c-1}+k)-t_{2c}\}$, we obtain

\al{
&|\E[h_{c,k}(\bm{X}_{t_{1},k},\dots,\bm{X}_{t_{c},k})h_{c,k}(\bm{X}_{t_{c+1},k},\dots,\bm{X}_{t_{2c},k})]|\notag\\
\leq{}& 12M_c^{\frac{2}{2+\delta}}\beta^{\frac{\delta}{2+\delta}}(\min\{t_{2c}-t_{2c-1},(t_{2c-1}+k)-t_{2c}\}).
}

Combining these inequalities yields (iii), which concludes the proof.

\ei
\ei

\hfill $\Box$

\newpage

\subsection{Proofs of results from Section~\ref{sec4}}\label{sec:5.3}

\subsubsection{Proof of Lemma \ref{cum}} \label{sec:proofcum}
By Theorem 5.2 of \cite{bradley05}, $(V_t)_{t\in\Z}$ is $\alpha$-mixing with mixing coefficients
\als{\label{amixing}
&\alpha^{V}(m)\leq \alpha^{X}(m)+\alpha^{X^{(1)}}(m)+\dots+\alpha^{X^{(q)}}(m)=(q+1)\alpha^{X}(m),\notag
}
where the latter identity is due to the definition of the processes $(X_t^{(j)})_{t\in\Z}$, $j=1,\dots,q$. Following the arguments in the proof of Lemma 4.1 in \cite{kley14}, we obtain the result. \hfill $\Box$


\subsubsection{Proof of Lemma \ref{auxlemmavar}} \label{sec:proofauxlemmavar}
For notational convenience, let $\bm{X}_{t_j,k_j}:=(X_{t_j},X_{t_j+k_j})^T$. Note that,
\al{
\E\Big[\prod_{j=1}^qh_{1,k_j}(\bm{X}_{t_j,k_j})\Big]={}&\E\bigg[\prod_{j=1}^{q}\E\Big[h\Big(\bm{X}_{t_j,k_j},\bm{X}_{t_{j},k_j}^{(1,j)},\dots,\bm{X}_{t_{j},k_j}^{(m-1,j)}\Big)-\xi_{k_j}\Big|\bm{X}_{t_j,k_j}\Big]\Big].\notag
}
Define 
\[
\mathcal{G}_j:=\sigma\Big(\bm{X}_{t_1,k_1},\dots,\bm{X}_{t_q,k_q},\bm{X}_{t_{j+1},k_{j+1}}^{(1,j+1)},\dots,\bm{X}_{t_{q},k_q}^{(1,q)},\dots,\bm{X}_{t_{j+1},k_{j+1}}^{(m-1,j+1)},\dots,\bm{X}_{t_{q},k_q}^{(m-1,q)}\Big).
\]
By the law of the total expectation we obtain
\al{
&\E\bigg[\prod_{j=1}^{q}\Big(h\Big(\bm{X}_{t_j,k_j},\bm{X}_{t_{j},k_j}^{(1,j)},\dots,\bm{X}_{t_{j},k_j}^{(m-1,j)}\Big)-\xi_{k_j}\Big)\Big]\\
&\hspace{1.5cm}=\E\bigg[\E\Big[\Big(h\Big(\bm{X}_{t_1,k_1},\bm{X}_{t_{1},k_1}^{(1,1)},\dots,\bm{X}_{t_{1},k_1}^{(m-1,1)}\Big)-\xi_{k_1}\Big)\\
&\hspace{3cm}\prod_{j=2}^{q}\Big(h\Big(\bm{X}_{t_j,k_j},\bm{X}_{t_{j},k_j}^{(1,j)},\dots,\bm{X}_{t_{j},k_j}^{(m-1,j)}\Big)-\xi_{k_j}\Big)\Big|\mathcal{G}_1\Big]\Big]\\
&\hspace{1.5cm}=\E\bigg[\E\Big[\Big(h\Big(\bm{X}_{t_1,k_1},\bm{X}_{t_{1},k_1}^{(1,1)},\dots,\bm{X}_{t_{1},k_1}^{(m-1,1)}\Big)-\xi_{k_1}\Big)\Big|\mathcal{G}_1\Big]\\
&\hspace{3cm}\prod_{j=2}^{q}\Big(h\Big(\bm{X}_{t_j,k_j},\bm{X}_{t_{j},k_j}^{(1,j)},\dots,\bm{X}_{t_{j},k_j}^{(m-1,j)}\Big)-\xi_{k_j}\Big)\Big],
}
where the latter inequality follows as $$\prod_{j=2}^{q}\big(h(\bm{X}_{t_j,k_j},\bm{X}_{t_{j},k_j}^{(1,j)},\dots,\bm{X}_{t_{j},k_j}^{(m-1,j)})-\xi_{k_j}\big)$$ is $\mathcal{G}_1$-measurable. Moreover, $\bm{X}_{t_1,k_1}$ is $\mathcal{G}_1$-measurable and $\sigma(\bm{X}_{t_{1},k_1}^{(1,1)},\dots,\bm{X}_{t_{1},k_1}^{(m-1,1)})$ is independent of $\mathcal{G}_1$. From the property of the conditional expectation stated in the proof of Lemma \ref{lemma1m}, it follows that
\al{
&\E\Big[\Big(h\Big(\bm{X}_{t_1,k_1},\bm{X}_{t_{1},k_1}^{(1,1)},\dots,\bm{X}_{t_{1},k_1}^{(m-1,1)}\Big)-\xi_{k_1}\Big)\Big|\mathcal{G}_1\Big]\\
={}&\E\Big[\Big(h\Big(\bm{X}_{t_1,k_1},\bm{X}_{t_{1},k_1}^{(1,1)},\dots,\bm{X}_{t_{1},k_1}^{(m-1,1)}\Big)-\xi_{k_1}\Big)\Big|\sigma (\bm{X}_{t_1,k_1})\Big]
={}h_{1,k_1}(\bm{X}_{t_1,k_1}).
}
Hence, the same arguments as above yield
\al{
&\E\bigg[\prod_{j=1}^{q}\Big(h\Big(\bm{X}_{t_j,k_j},\bm{X}_{t_{j},k_j}^{(1,j)},\dots,\bm{X}_{t_{j},k_j}^{(m-1,j)}\Big)-\xi_{k_j}\Big)\Big]\\
={}&\E\bigg[h_{1,k_1}(\bm{X}_{t_1,k_1})\prod_{j=2}^{q}\Big(h\Big(\bm{X}_{t_j,k_j},\bm{X}_{t_{j},k_j}^{(1,j)},\dots,\bm{X}_{t_{j},k_j}^{(m-1,j)}\Big)-\xi_{k_j}\Big)\Big]\\
={}&\E\bigg[\E\Big[\Big(h\Big(\bm{X}_{t_2,k_2},\bm{X}_{t_{2},k_2}^{(2,2)},\dots,\bm{X}_{t_{2},k_2}^{(m-1,2)}\Big)-\xi_{k_2}\Big)\Big|\mathcal{G}_2\Big]\\
&\qquad\qquad \cdot h_{1,k_1}(\bm{X}_{t_1,k_1})\prod_{j=3}^{q}\Big(h\Big(\bm{X}_{t_j,k_j},\bm{X}_{t_{j},k_j}^{(1,j)},\dots,\bm{X}_{t_{j},k_j}^{(m-1,j)}\Big)-\xi_{k_j}\Big)\Big]\\
={}&\E\bigg[h_{1,k_2}(\bm{X}_{t_2,k_2})h_{1,k_1}(\bm{X}_{t_1,k_1})\prod_{j=3}^{q}\Big(h\Big(\bm{X}_{t_j,k_j},\bm{X}_{t_{j},k_j}^{(1,j)},\dots,\bm{X}_{t_{j},k_j}^{(m-1,j)}\Big)-\xi_{k_j}\Big)\Big].
}

Repeating these steps $q-1$ times yields the result.\\

\bi
\item[(i)] Under (C0) we have $h_{1,k}^{\tau}\colvec{x_1}{y_1}=\E\Big[h^{\tau}\Big(\colvec{x_1}{y_1},\colvec{X_0}{X_k}\Big)\Big] $ with
\al{
h^{\tau}\Big(\colvec{x_1}{y_1},\colvec{x_2}{y_2}\Big)=4\Big(I(x_1<x_2)-\frac{1}{2}\Big)\Big(I(y_1<y_2)-\frac{1}{2}\Big).
}
Note that $\E[Y_t^{(j)}]=0$ under (C0), we obtain from (\ref{eprodhm}) that
\als{\label{covcum}
&\Cov\Big(h_{1,k_1}^{\tau}\colvec{X_{t_1}}{X_{t_1+k_1}},h_{1,k_2}^{\tau}\colvec{X_{t_2}}{X_{t_2+k_2}}\Big)\notag\\
={}&\E\Big[\Big(h^{\tau}\Big(\colvec{X_{t_1}}{X_{t_1+k_1}},\colvec{X_{t_1}^{(1)}}{X_{t_1+k_1}^{(1)}}\Big)-\tau_{k_1}\Big)\Big(h^{\tau}\Big(\colvec{X_{t_2}}{X_{t_2+k_2}},\colvec{X_{t_2}^{(2)}}{X_{t_2+k_2}^{(2)}}\Big)-\tau_{k_2}\Big)\Big]\notag\\
={}&\E\Big[\Big(4\Big(I(X_{t_1}<X_{t_1}^{(1)})-\frac{1}{2}\Big)\Big(I(X_{t_1+k_1}<X_{t_1+k_1}^{(1)})-\frac{1}{2}\Big)-\tau_{k_1}\Big)\notag\\
&\cdot\Big(4\Big(I(X_{t_2}<X_{t_2}^{(2)})-\frac{1}{2}\Big)\Big(I(X_{t_2+k_2}<X_{t_2+k_2}^{(2)})-\frac{1}{2}\Big)-\tau_{k_2}\Big)\Big]\notag\\
={}&\E[(4\yt_{t_1}\yt_{t_1+k_1}-\tau_{k_1})(4\ytt_{t_2}\ytt_{t_2+k_2}-\tau_{k_2})]\notag\\
={}&16\E[\yt_{t_1}\yt_{t_1+k_1}\ytt_{t_2}\ytt_{t_2+k_2}]-4\E[\yt_{t_1}\yt_{t_1+k_1}]\tau_{k_2}-4\E[\ytt_{t_2}\ytt_{t_2+k_2}]\tau_{k_1}+\tau_{k_1}\tau_{k_2}.\notag\\
}
The equivalent representation of moments in terms of cumulants yields
\al{
\E[\yt_{t_1}&\yt_{t_1+k_1}\ytt_{t_2}\ytt_{t_2+k_2}]=\cum (\yt_{t_1},\yt_{t_1+k_1},\ytt_{t_2},\ytt_{t_2+k_2}) \\
&+\cum (\yt_{t_1},\yt_{t_1+k_1})\cum(\ytt_{t_2},\ytt_{t_2+k_2})+\cum(\yt_{t_1},\ytt_{t_2})\cum(\yt_{t_1+k_1},\ytt_{t_2+k_2})\\
&+\cum(\yt_{t_1},\ytt_{t_2+k_2})\cum(\yt_{t_1+k_1},\ytt_{t_2})
}
For all $t,k\in\Z$ and $l=1,2$, $\E[Y_t^{(l)}]=0$ and we have
\als{\label{cumtau}
&\cum (Y_{t}^{(l)},Y_{t+k}^{(l)})=\E[Y_{t}^{(l)}Y_{t+k}^{(l)}]\notag\\
={}&\int_{\R^2}\Big(I(x_{t}<x_{t}^{(l)})-\frac{1}{2}\Big)\Big(I(x_{t+k}<x_{t+k}^{(l)})-\frac{1}{2}\Big)dF_k\colvec{x_t^{(l)}}{x_{t+k}^{(l)}}dF_k\colvec{x_t}{x_{t+k}}\notag\\
={}&\int_{\R^2}F_k\colvec{x_t}{x_{t+k}}dF_k\colvec{x_t}{x_{t+k}}-\frac{1}{4}\notag\\
={}&\int_{[0,1]^2}C_k(u,v)dC_k(u,v)-\frac{1}{4} 
=\frac{1}{4}\tau_k
}
and
\als{\label{cumrho}
&\cum (Y_{t}^{(1)},Y_{t+k}^{(2)})=\E[Y_{t}^{(1)}Y_{t+k}^{(2)}]\notag\\
={}&\int_{\R^2}\Big(I(x_{t}<x_{t}^{(1)})-\frac{1}{2}\Big)\Big(I(x_{t+k}<x_{t+k}^{(2)})-\frac{1}{2}\Big)dF(x_t^{(1)})dF(x_{t+k}^{(2)})dF_k\colvec{x_t}{x_{t+k}}\notag\\
={}&\int_{\R^2}F(x_t)F(x_{t+k})dF_k\colvec{x_t}{x_{t+k}}-\frac{1}{4}\notag\\
={}&\int_{[0,1]^2}uvdC_k(u,v)-\frac{1}{4} = \frac{1}{12}\rho(k)
}
where $C_k$ is the copula associated with $(X_t,X_{t+k})$ [see e.g. \cite{schmid10}] and $\rho(k)$ is the population version of Spearman's $\rho$ at lag $k$. Hence,
\als{\label{4thmoment}
\E[\yt_{t_1}&\yt_{t_1+k_1}\ytt_{t_2}\ytt_{t_2+k_2}]=\cum (\yt_{t_1},\yt_{t_1+k_1},\ytt_{t_2},\ytt_{t_2+k_2})+\frac{1}{16}\tau_{k_1}\tau_{k_2}\notag\\
&+\frac{1}{144}\rho(t_2-t_1)\rho(t_2+k_2-(t_1+k_1))+\frac{1}{144}\rho(t_2+k_2-t_1)\rho(t_2-(t_1+k_1))\notag\\
{}
}
and inserting equations (\ref{cumtau}) and (\ref{4thmoment}) in (\ref{covcum}) yield the result.

\item[(ii)]

Under (C0) we have $h_{1,k}^{\rho}\colvec{x}{y}=\E\Big[h^{\rho}\Big(\colvec{x}{y},\colvec{X_0^{(1)}}{X_k^{(1)}},\colvec{X_0^{(2)}}{X_k^{(2)}}\Big)\Big]$ with
 \al{
&h^{\rho}\Big(\colvec{x_1}{y_1},\colvec{x_2}{y_2},\colvec{x_3}{y_3}\Big)={}\sum_{\gamma\in\Gamma\{1,2,3\}}2\Big(I(x_{\gamma(1)}<x_{\gamma(2)})-\frac{1}{2}\Big)\Big(I(y_{\gamma(1)}<y_{\gamma(3)})-\frac{1}{2}\Big).
}
The first order kernel being centered by definition of the Hoeffding decomposition, from (\ref{eprodhm}) we know that 
\al{
\Cov\Big(h_{1,k_1}^{\rho}&\colvec{X_{t_1}}{X_{t_1+k_1}},h_{1,k_2}^{\rho}\colvec{X_{t_2}}{X_{t_2+k_2}}\Big)\\
={}&\E\Big[\Big(h^{\rho}\Big(\colvec{X_{t_1}^{(1)}}{X_{t_1+k_1}^{(1)}},\colvec{X_{t_1}^{(2)}}{X_{t_1+k_1}^{(2)}},\colvec{X_{t_1}^{(3)}}{X_{t_1+k_1}^{(3)}}\Big)-\rho_{k_1}\Big)\\
&\qquad\cdot\Big(h^{\rho}\Big(\colvec{X_{t_2}^{(1)}}{X_{t_2+2}^{(1)}},\colvec{X_{t_2}^{(4)}}{X_{t_2+k_2}^{(4)}},\colvec{X_{t_2}^{(5)}}{X_{t_2+k_2}^{(5)}}\Big)-\rho_{k_2}\Big)\Big].
}
Thus, as $\E[I(X_{t}^{(i)}<X_{t}^{(j)})-\frac{1}{2}]=0$ under (C0) for any $i,j=1,\dots,5;\,i\neq j$, we obtain from (\ref{eprodhm}) that
\als{\label{covdecrho}
&\Cov\Big(h_{1,k_1}^{\rho}\colvec{X_{t_1}}{X_{t_1+k_1}},h_{1,k_2}^{\rho}\colvec{X_{t_2}}{X_{t_2+k_2}}\Big)\notag\\
={}&4\sum_{\gamma\in\Gamma\{1,2,3\}}\sum_{\gt\in\Gamma\{1,4,5\}}\Big\{\E\Big(I(X_{t_1}^{(\gamma(1))}<X_{t_1}^{(\gamma(2))})-\frac{1}{2}\Big)\Big(I(X_{t_1+k_1}^{(\gamma(1))}<X_{t_1+k_1}^{(\gamma(3))})-\frac{1}{2}\Big)\notag\\
&\qquad\cdot\Big(I(X_{t_2}^{(\gt(1))}<X_{t_2}^{(\gt(2))})-\frac{1}{2}\Big)\Big(I(X_{t_2+k_2}^{(\gt(1))}<X_{t_2+k_2}^{(\gt(3))})-\frac{1}{2}\Big)\Big\}-\rho_{k_1}\rho_{k_2}\notag\\
={}&4\sum_{\gamma\in\Gamma\{1,2,3\}}\sum_{\gt\in\Gamma\{1,4,5\}}\Big[\cum\Big(I(X_{t_1}^{(\gamma(1))}<X_{t_1}^{(\gamma(2))}),I(X_{t_1+k_1}^{(\gamma(1))}<X_{t_1+k_1}^{(\gamma(3))}),\notag\\
&\hspace{4cm} I(X_{t_2}^{(\gt(1))}<X_{t_2}^{(\gt(2))}),I(X_{t_2+k_2}^{(\gt(1))}<X_{t_2+k_2}^{(\gt(3))})\Big)\notag\\
&\hspace{1cm}+\cum\Big(I(X_{t_1}^{(\gamma(1))}<X_{t_1}^{(\gamma(2))})-\frac{1}{2},I(X_{t_2}^{(\gt(1))}<X_{t_2}^{(\gt(2))})-\frac{1}{2}\Big)\notag\\
&\hspace{1cm}\qquad\cdot \cum\Big(I(X_{t_1+k_1}^{(\gamma(1))}<X_{t_1+k_1}^{(\gamma(3))})-\frac{1}{2},I(X_{t_2+k_2}^{(\gt(1))}<X_{t_2+k_2}^{(\gt(3))})-\frac{1}{2}\Big)\notag\\
&\hspace{1cm}+\cum\Big(I(X_{t_1}^{(\gamma(1))}<X_{t_1}^{(\gamma(2))})-\frac{1}{2},I(X_{t_2+k_2}^{(\gt(1))}<X_{t_2+k_2}^{(\gt(3))})-\frac{1}{2}\Big)\notag\\
&\hspace{1cm}\qquad\cdot\cum\Big(I(X_{t_1+k_1}^{(\gamma(1))}<X_{t_1+k_1}^{(\gamma(3))})-\frac{1}{2},I(X_{t_2}^{(\gt(1))}<X_{t_2}^{(\gt(2))})-\frac{1}{2}\Big)\Big],
}
where we have used the representation of centered fourth moments in terms of cumulants, property (v) of Theorem 2.3.1 in \cite{brillinger75} and (\ref{cumrho}) 
\al{
&\cum\Big(I(X_{t_1}^{(\gamma(1))}<X_{t_1}^{(\gamma(2))})-\frac{1}{2},I(X_{t_1+k_1}^{(\gamma(1))}<X_{t_1+k_1}^{(\gamma(3))})-\frac{1}{2}\Big)\\
&\qquad\cdot\cum\Big(I(X_{t_2}^{(\gt(1))}<X_{t_2}^{(\gt(2))})-\frac{1}{2},I(X_{t_2+k_2}^{(\gt(1))}<X_{t_2+k_2}^{(\gt(3))})-\frac{1}{2}\Big)=\frac{1}{144}\rho_{k_1}\rho_{k_2}.
}
Furthermore, the only permutations $\gamma$ and $\gt$ for which not all products of second order cumulants in (\ref{covdecrho}) contain one cumulant with one independent factor and thus equal $0$ are those with $\gamma(1)=\gt(1)=1$. For each of these $4$ combinations we obtain
\als{\label{cum1}
&\cum\Big(I(X_{t_1}^{(\gamma(1))}<X_{t_1}^{(\gamma(2))})-\frac{1}{2},I(X_{t_2}^{(\gt(1))}<X_{t_2}^{(\gt(2))})-\frac{1}{2}\Big)\notag\\
&\hspace{1cm}\qquad\cdot \cum\Big(I(X_{t_1+k_1}^{(\gamma(1))}<X_{t_1+k_1}^{(\gamma(3))})-\frac{1}{2},I(X_{t_2+k_2}^{(\gt(1))}<X_{t_2+k_2}^{(\gt(3))})-\frac{1}{2}\Big)\notag\\
={}&\frac{1}{144}\rho(t_2-t_1)\rho(t_2+k_2-(t_1+k_1))
}
and
\als{\label{cum2}
&\cum\Big(I(X_{t_1}^{(\gamma(1))}<X_{t_1}^{(\gamma(2))})-\frac{1}{2},I(X_{t_2+k_2}^{(\gt(1))}<X_{t_2+k_2}^{(\gt(3))})-\frac{1}{2}\Big)\notag\\
&\hspace{1cm}\qquad\cdot\cum\Big(I(X_{t_1+k_1}^{(\gamma(1))}<X_{t_1+k_1}^{(\gamma(3))})-\frac{1}{2},I(X_{t_2}^{(\gt(1))}<X_{t_2}^{(\gt(2))})-\frac{1}{2}\Big)\notag\\
={}&\frac{1}{144}\rho(t_2+k_2-t_1)\rho(t_2-(t_1+k_1)).
}

Plugging (\ref{cum1}) and (\ref{cum2}) into (\ref{covdecrho}) concludes the proof.
\ei

\hfill $\Box$

\newpage

\subsubsection{Proof of \eqref{laggedyoshihara}} \label{sec:proofyoshi}
We will prove this result only for positive lags $k$ as the proof for negative lags is analogous. More precisely we consider
\als{\label{covariancehck}
&\E\Big[\Big(\binom{n-k}{c}U_{n-k}^{(c)}\Big)^2\Big]\notag\\
={}&\sum\limits_{1\leq t_1<\dots<t_c\leq n-k}\sum\limits_{1\leq t_{c+1}<\dots<t_{2c}\leq n-k}\E[h_{c,k}(\bm{X}_{t_{1},k},\dots,\bm{X}_{t_{c},k})h_{c,k}(\bm{X}_{t_{c+1},k},\dots,\bm{X}_{t_{2c},k})]
}
and prove that
\bea
\sup_{0\leq k\leq\floor{r_n}}\E\Big[\Big(\binom{n-k}{c}U_{n-k}^{(c)}\Big)^2\Big]=O(n^{2c-1-\theta})
\eea

Finally, using that $\inf_{0\leq k\leq\floor{r_n}}\binom{n-k}{c}\geq Kn^{c}$ for some constant $K$, establishes (\ref{laggedyoshihara}).\\

For any fixed $0\leq k\leq\floor{r_n}$, decompose (\ref{covariancehck}) into sums according to the following 3 cases:

\begin{itemize}
\item[(1)] all $2c$ indices are different,
\item[(2)] $2c-1$ indices are different or
\item[(3)] $2(c-1)$ or less indices are different,
\ei
that is
\bea
&&\Big|\sum\limits_{1\leq t_1<\dots<t_c\leq n-k}\sum\limits_{1\leq t_{c+1}<\dots<t_{2c}\leq n-k}\E[h_{c,k}(\bm{X}_{t_{1},k},\dots,\bm{X}_{t_{c},k})h_{c,k}(\bm{X}_{t_{c+1},k},\dots,\bm{X}_{t_{2c},k})]\Big|\\
&\leq&\sum\limits_{\substack{1\leq t_1<\dots<t_c\leq n-k,\,1\leq t_{c+1}<\dots<t_{2c}\leq n-k\\ \text{case (1)}}}|\E[h_{c,k}(\bm{X}_{t_{1},k},\dots,\bm{X}_{t_{c},k})h_{c,k}(\bm{X}_{t_{c+1},k},\dots,\bm{X}_{t_{2c},k})]|\\
&&+\sum\limits_{\substack{1\leq t_1<\dots<t_c\leq n-k,\,1\leq t_{c+1}<\dots<t_{2c}\leq n-k\\ \text{case (2)}}}|\E[h_{c,k}(\bm{X}_{t_{1},k},\dots,\bm{X}_{t_{c},k})h_{c,k}(\bm{X}_{t_{c+1},k},\dots,\bm{X}_{t_{2c},k})]|\\
&&+\sum\limits_{\substack{1\leq t_1<\dots<t_c\leq n-k,\,1\leq t_{c+1}<\dots<t_{2c}\leq n-k\\ \text{case (3)}}}|\E[h_{c,k}(\bm{X}_{t_{1},k},\dots,\bm{X}_{t_{c},k})h_{c,k}(\bm{X}_{t_{c+1},k},\dots,\bm{X}_{t_{2c},k})]|
\eea
In the sequel, denote by $t_{(j)}$ the $j$-th smallest of all distinct indices among $t_1,\dots,t_{2c}$.\\

In \textbf{case (1)} we distinguish the following cases:

\bi
\item[(1.1)] $(t_{(2)}-t_{(1)})>k$ or $(t_{(2c)}-t_{(2c-1)})>k$.
\item[(1.2)] $(t_{(2)}-t_{(1)})\leq k$ and $(t_{(2c)}-t_{(2c-1)})\leq k$.
\ei

In \textbf{case (1.1)}, consider the set
\al{
\mathcal{S}_k^{(1.1)}(v):={}&\{t_1,\dots,t_{2c}:1\leq t_1 <\dots <t_c;1\leq t_{c+1} <\dots <t_{2c};
 t_1\neq \dots \neq t_{2c};\\
&(t_{(2)}-t_{(1)})>k\text{ or }(t_{(2c)}-t_{(2c-1)})>k;\max\{t_{(2)}-t_{(1)},t_{(2c)}-t_{(2c-1)}\}=v\}
}
and observe that $\#\mathcal{S}_{k}^{(1.1)}(v)\leq  (v+k)n^{2(c-1)}$, where $\#S$ denotes the cardinality of the set $S$.\\

Then, we obtain from Lemma \ref{boundshm} (2) (ii), for some constant $K$,
\als{\label{case(1.1)}
&\sum\limits_{\substack{1\leq t_1<\dots<t_c\leq n-k,\,1\leq t_{c+1}<\dots<t_{2c}\leq n-k\\ t_1\neq\cdots\neq t_{2c}\\ (t_{(2)}-t_{(1)})>k\text{ or }(t_{(2c)}-t_{(2c-1)})>k}}|\E[h_{c,k}(\bm{X}_{t_{1},k},\dots,\bm{X}_{t_{c},k})h_{c,k}(\bm{X}_{t_{c+1},k},\dots,\bm{X}_{t_{2c},k})]|\notag\\
\leq{}&K\sum_{v=1}^{n-k}\sum_{t_1,\dots,t_{2c}\in\mathcal{S}_{k}^{(1.1)}(v)}\beta^{\frac{\delta}{2+\delta}}(\max\{t_{(2)}-(t_{(1)}+k),t_{(2c)}-(t_{(2c-1)}+k)\})\notag\\
\leq{}&K\sum_{v=1}^{n-k}\beta^{\frac{\delta}{2+\delta}}(v)\#\mathcal{S}_{k}^{(1.1)}(v)
\leq{}Kn^{2(c-1)}\sum_{v=1}^{n-k}(v+k)\beta^{\frac{\delta}{2+\delta}}(v)\notag\\
\leq{}&O(n^{1-\theta})n^{2(c-1)}
={}O(n^{2c-1-\theta}),
}
where we have bounded $\sum_{v=1}^{n}v\beta^{\frac{\delta}{2+\delta}}(v)$ from above by an integral and then concluded with Assumption (C3).\\ 

Next, in \textbf{case (1.2)}, 
\al{
&\sum\limits_{\substack{1\leq t_1<\dots<t_c\leq n-k,\,1\leq t_{c+1}<\dots<t_{2c}\leq n-k\\ t_1\neq\cdots\neq t_{2c}\\ (t_{(2)}-t_{(1)})\leq k\text{ and }(t_{(2c)}-t_{(2c-1)})\leq k}}|\E[h_{c,k}(\bm{X}_{t_{1},k},\dots,\bm{X}_{t_{c},k})h_{c,k}(\bm{X}_{t_{c+1},k},\dots,\bm{X}_{t_{2c},k})]|\\
={}&\sum\limits_{\substack{1\leq t_1<\dots<t_c\leq n-k,\,1\leq t_{c+1}<\dots<t_{2c}\leq n-k\\ t_1\neq\cdots\neq t_{2c}\\ (t_{(2)}-t_{(1)})\leq k\text{ and }(t_{(2c)}-t_{(2c-1)})\leq k\\
\text{and } ((t_{(3)}-t_{(2)})\leq 2k\text{ or }(t_{(2c-1)}-t_{(2(c-1))})\leq 2k)}}|\E[h_{c,k}(\bm{X}_{t_{1},k},\dots,\bm{X}_{t_{c},k})h_{c,k}(\bm{X}_{t_{c+1},k},\dots,\bm{X}_{t_{2c},k})]|\\
&+\sum\limits_{\substack{1\leq t_1<\dots<t_c\leq n-k,\,1\leq t_{c+1}<\dots<t_{2c}\leq n-k\\ t_1\neq\cdots\neq t_{2c}\\ (t_{(2)}-t_{(1)})\leq k\text{ and }(t_{(2c)}-t_{(2c-1)})\leq k\\
\text{and } (t_{(3)}-t_{(2)})>2k\text{ and }(t_{(2c-1)}-t_{(2(c-1))})>2k}}|\E[h_{c,k}(\bm{X}_{t_{1},k},\dots,\bm{X}_{t_{c},k})h_{c,k}(\bm{X}_{t_{c+1},k},\dots,\bm{X}_{t_{2c},k})]|\\
={}&(I)+(II).
}

For (I), we apply similar arguments as in the proof of Lemma \ref{boundshm} (2) (iii). That is, if we replace one after another all random variables by independent copies and then rereplace them by independent pairs with cdf $F_k$. We have for any permutation $\gamma$ of $\{1,\dots,2c\}$,
\al{
&|\E[h_{c,k}(\bm{X}_{t_{\gamma(1)},k},\dots,\bm{X}_{t_{\gamma(c)},k})h_{c,k}(\bm{X}_{t_{\gamma(c+1)},k},\dots,\bm{X}_{t_{\gamma(2c)},k})]|\\
\leq{}& K_cM_c^{\frac{2}{2+\delta}}\beta^{\frac{\delta}{2+\delta}}(\min_{\substack{i,j=1,\dots,2c\\i\neq j}}\{|t_{(j)}-t_{(i)}|,|(t_{(j)}+k)-t_{(i)}|\}).
}
for a constant $K_c$ depening only on $c$. Next, let $u(t_{(1)},\dots,t_{(2c)}):=\min_{\substack{i,j=1,\dots,2c\\i\neq j}}\{|t_{(j)}-t_{(i)}|,|(t_{(j)}+k)-t_{(i)}|\}$ which is always smaller than $k$ and consider the set
\al{
\mathcal{S}_k^{(I)}(v):={}&\{t_1,\dots,t_{2c}:1\leq t_1 <\dots <t_c;1\leq t_{c+1} <\dots <t_{2c};
 t_1\neq \dots \neq t_{2c};\\
&(t_{(2)}-t_{(1)})\leq k\text{ and }(t_{(2c)}-t_{(2c-1)})\leq k;(t_{(3)}-t_{(2)})\leq 2k\\
&\text{ or }(t_{(2c-1)}-t_{(2c-2)})\leq 2k; u(t_{(1)},\dots,t_{(2c)})=v\}.
}
Then, for some constant $K$,
\al{
|(I)|\leq{}& K\sum_{v=0}^{k}\sum_{t_1,\dots, t_{2c}\in\mathcal{S}_k^{(I)}(v)}\beta^{\frac{\delta}{2+\delta}}(u(t_{(1)},\dots,t_{(2c)}))\\
\leq{}&K\sum_{v=0}^{k}\beta^{\frac{\delta}{2+\delta}}(v)\#\mathcal{S}_k^{(I)}(v) 
\leq{} O(1)r_n^{2}n^{2c-3}
}
since $\sup\limits_{v=0,\dots,k}\#\mathcal{S}_k^{(I)}(v)\leq 2r_n^{2}n^{2c-3}$. Hence,
$
(I)=O(r_n^2n^{2c-3})=o(n^{2c-1-\theta}).
$
From Lemma \ref{boundshm} (2) (iii) we know that 
\al{
(II)\leq{}& 12M_c^{\frac{2}{2+\delta}}\sum\limits_{\substack{1\leq t_1<\dots<t_c\leq n-k,\,1\leq t_{c+1}<\dots<t_{2c}\leq n-k\\ t_1\neq\cdots\neq t_{2c}\\ (t_{(2)}-t_{(1)})\leq k\text{ and }(t_{(2c)}-t_{(2c-1)})\leq k\\ \text{and } (t_{(3)}-t_{(2)})>2k\text{ and }(t_{(2c-1)}-t_{(2(c-1))})>2k}}\beta^{\frac{\delta}{
2+\delta}}(m(t_1,\dots,t_{2c})).
} 
Consider the set
\al{
\mathcal{S}_k^{(II)}(v):={}&\{t_1,\dots,t_{2c}:1\leq t_1 <\dots <t_c;1\leq t_{c+1} <\dots <t_{2c};
 t_1\neq \dots \neq t_{2c};\\
&(t_{(2)}-t_{(1)})\leq k\text{ and }(t_{(2c)}-t_{(2c-1)})\leq k;(t_{(3)}-t_{(2)})> 2k\\
&\text{ and }(t_{(2c-1)}-t_{(2c-2)})> 2k;
m(t_{(1)},\dots,t_{(2c)})=v\}
}
with $\#\mathcal{S}_{k}^{(II)}(v)\leq  2(v+1)n^{2c-2}$. Then, for some constant $K$,

\al{
|(II)|\leq{}& K\sum_{v=0}^{k}\sum_{t_1,\dots, t_{2c}\in\mathcal{S}_k^{(II)}(v)}\beta^{\frac{\delta}{2+\delta}}(m(t_{(1)},\dots,t_{(2c)}))
\leq{}K\sum_{v=0}^{k}\beta^{\frac{\delta}{2+\delta}}(v)\#\mathcal{S}_k^{(II)}(v)\\
\leq{}&K\sum_{v=0}^{k}2(v+1)\beta^{\frac{\delta}{2+\delta}}(v)n^{2c-2}
\leq{}O(r_n^{1-\theta})n^{2c-2}
={}O(r_n^{1-\theta}n^{2c-2})
}
and hence,
$
(II)=O(r_n^{1-\theta}n^{2c-2})=O(n^{2c-1-\theta}).
$

Therefore, in case (1.2) we have
\als{\label{case(1.2)}
\sum\limits_{\substack{1\leq t_1<\dots<t_c\leq n-k,\,1\leq t_{c+1}<\dots<t_{2c}\leq n-k\\ t_1\neq\cdots\neq t_{2c}\\ (t_{(2)}-t_{(1)})\leq k\text{ and }(t_{(2c)}-t_{(2c-1)})\leq k}}|\E[h_{c,k}(\bm{X}_{t_{1},k},\dots,\bm{X}_{t_{c},k})&h_{c,k}(\bm{X}_{t_{c+1},k},\dots,\bm{X}_{t_{2c},k})]|\notag\\
={}&O(n^{2c-1-\theta}).
}
Combining equations (\ref{case(1.1)}) and (\ref{case(1.2)}) yields
\als{\label{case1m}
\sum\limits_{\substack{1\leq t_1<\dots<t_c\leq n-k,\,1\leq t_{c+1}<\dots<t_{2c}\leq n-k\\ \text{case (1)}}}|\E[h_{c,k}(\bm{X}_{t_{1},k},\dots,\bm{X}_{t_{c},k})h_{c,k}(\bm{X}_{t_{c+1},k},\dots,\bm{X}_{t_{2c},k})]|=O(n^{2c-1-\theta}),
}
which concludes the consideration of case (1).\\

In \textbf{case (2)}, we encounter the following situations:

\bi
\item[(2.1)] the index appearing twice is not $t_{(1)}$ or $t_{(2)}$.
\item[(2.2)] the index appearing twice is $t_{(1)}$ or $t_{(2)}$. 
\ei

Then,
\al{
&\sum\limits_{\substack{1\leq t_1<\dots<t_c\leq n-k,\,1\leq t_{c+1}<\dots<t_{2c}\leq n-k\\ \text{case (2)}}}|\E[h_{c,k}(\bm{X}_{t_{1},k},\dots,\bm{X}_{t_{c},k})h_{c,k}(\bm{X}_{t_{c+1},k},\dots,\bm{X}_{t_{2c},k})]|\\
&=\Big(\sum\limits_{\substack{1\leq t_1<\dots<t_c\leq n-k,\,1\leq t_{c+1}<\dots<t_{2c}\leq n-k\\ \text{case (2.1)}}}+\sum\limits_{\substack{1\leq t_1<\dots<t_c\leq n-k,\,1\leq t_{c+1}<\dots<t_{2c}\leq n-k\\ \text{case (2.2)}}}\\
&\hspace{2cm}|\E[h_{c,k}(\bm{X}_{t_{1},k},\dots,\bm{X}_{t_{c},k})h_{c,k}(\bm{X}_{t_{c+1},k},\dots,\bm{X}_{t_{2c},k})]|\\
&=:{}S_1+S_2
}
In \textbf{case (2.1)}, consider the following situations:
\bi
\item[(a)] $t_{(2)}-t_{(1)}>k$
\item[(b)] $t_{(2)}-t_{(1)}\leq k$ and $t_{(3)}-t_{(2)}>2k$
\item[(c)] $t_{(2)}-t_{(1)}\leq k$ and $t_{(3)}-t_{(2)}\leq 2k$
\ei
In \textbf{situation (a)}, similarly as in the proof of Lemma \ref{boundshm} (ii), we replace the pair with the smallest index $\bm{X}_{t_{(1)},k}$ by an independent copy in order to bound the summand from above by $2M^{\frac{2}{2+\delta}}\beta^{\frac{\delta}{2+\delta}}(t_{(2)}-(t_{(1)}+k))$. Hence, by assumption (C3),
\al{
&\sum\limits_{\substack{1\leq t_1<\dots<t_c\leq n-k,\,1\leq t_{c+1}<\dots<t_{2c}\leq n-k\\ \text{case (2.1)}\text{ (a)} }}|\E[h_{c,k}(\bm{X}_{t_{1},k},\dots,\bm{X}_{t_{c},k})h_{c,k}(\bm{X}_{t_{c+1},k},\dots,\bm{X}_{t_{2c},k})]|\\
&\hspace{0cm}\leq n^{2c-3}2M_c^{\frac{2}{2+\delta}}\sum_{t_{(2)}-t_{(1)}>k}\beta^{\frac{\delta}{2+\delta}}(t_{(2)}-(t_{(1)}+k))\\
&\hspace{0cm}\leq n^{2c-2}2M_c^{\frac{2}{2+\delta}}\sum_{u=1}^{\infty}\beta^{\frac{\delta}{2+\delta}}(u)=O(n^{2(c-1)}).
}
Next, in \textbf{situation (b)}, with similar arguments as in the proof of Lemma \ref{boundshm} (iii), we replace one after another $X_{t_{(1)}}$, $X_{t_{(2)}}$, $X_{t_{(1)}+k}$ and $X_{t_{(2)}+k}$ by independent copies and obtain with assumption (C3),
\al{
&\sum\limits_{\substack{1\leq t_1<\dots<t_c\leq n-k,\,1\leq t_{c+1}<\dots<t_{2c}\leq n-k\\ \text{case (2.1)}\text{ (a)} }}|\E[h_{c,k}(\bm{X}_{t_{1},k},\dots,\bm{X}_{t_{c},k})h_{c,k}(\bm{X}_{t_{c+1},k},\dots,\bm{X}_{t_{2c},k})]|\\
&\hspace{0cm}\leq n^{2c-3}12M_c^{\frac{2}{2+\delta}}\sum_{t_{(2)}-t_{(1)}\leq k}\beta^{\frac{\delta}{2+\delta}}(\min\Big\{t_{(2)}-t_{(1)},t_{(1)}+k-t_{(2)}\Big\})\\
&\hspace{0cm}\leq n^{2c-2}12M_c^{\frac{2}{2+\delta}}\sum_{u=1}^{k}\beta^{\frac{\delta}{2+\delta}}(\min\Big\{u,k-u\Big\})\\
&\hspace{0cm}\leq n^{2c-2}12M_c^{\frac{2}{2+\delta}}2\sum_{u=0}^{\floor{\frac{k}{2}}}\beta^{\frac{\delta}{2+\delta}}(u)=O(n^{2(c-1)}).
}
In \textbf{situation (c)}, we use Lemma \ref{boundshm} (i) and the fact that in this case the number of summands is of order $O(r_n^2n^{2c-3})$, that is
\al{
&\sum\limits_{\substack{1\leq t_1<\dots<t_c\leq n-k,\,1\leq t_{c+1}<\dots<t_{2c}\leq n-k\\ \text{case (2.1)}\text{ (c)} }}|\E[h_{c,k}(\bm{X}_{t_{1},k},\dots,\bm{X}_{t_{c},k})h_{c,k}(\bm{X}_{t_{c+1},k},\dots,\bm{X}_{t_{2c},k})]|\\
&\hspace{3cm}\leq n^{2c-4}\sum_{\substack{t_{(2)}-t_{(1)}\leq k\\t_{(3)}-t_{(2)}\leq 2k}}M_c^{\frac{2}{2+\delta}}=O(r_n^2n^{2c-3})=O(n^{2(c-1)})
}
Therefore,
\al{
\sum\limits_{\substack{1\leq t_1<\dots<t_c\leq n-k,\,1\leq t_{c+1}<\dots<t_{2c}\leq n-k\\ \text{case (2.1)}}}|\E[h_{c,k}(\bm{X}_{t_{1},k},\dots,\bm{X}_{t_{c},k})h_{c,k}(\bm{X}_{t_{c+1},k},\dots,\bm{X}_{t_{2c},k})]|=O(n^{2(c-1)})
}
Since in \textbf{case (2.2)}, the index appearing twice is $t_{(1)}$ or $t_{(2)}$, the indices $t_{(2c-2)}$ and $t_{(2c-1)}$ appear only once. Thus the case can be handled by the similar arguments as case (2.1), i.e. we obtain 
\al{
\sum\limits_{\substack{1\leq t_1<\dots<t_c\leq n-k,\,1\leq t_{c+1}<\dots<t_{2c}\leq n-k\\ \text{case (2.2)}}}|\E[h_{c,k}(\bm{X}_{t_{1},k},\dots,\bm{X}_{t_{c},k})h_{c,k}(\bm{X}_{t_{c+1},k},\dots,\bm{X}_{t_{2c},k})]|=O(n^{2(c-1)}).
}
Cases (2.1) and (2.2) yield
\als{\label{case2m}
&\sum\limits_{\substack{1\leq t_1<\dots<t_c\leq n-k,\,1\leq t_{c+1}<\dots<t_{2c}\leq n-k\\ \text{case (2)}}}|\E[h_{c,k}(\bm{X}_{t_{1},k},\dots,\bm{X}_{t_{c},k})h_{c,k}(\bm{X}_{t_{c+1},k},\dots,\bm{X}_{t_{2c},k})]|\notag\\
&\hspace{4cm}=O(n^{2(c-1)})=O(n^{2c-1-\theta})
}
which concludes the consideration of case (2).\\

In \textbf{case (3)} observe that the number of summands is of order $O(n^{2(c-1)})$, such that  together with Lemma \ref{boundshm} (2) (i) we can conclude that 
\als{\label{case3m}
&\sum\limits_{\substack{1\leq t_1<\dots<t_c\leq n-k,\,1\leq t_{c+1}<\dots<t_{2c}\leq n-k\\ \text{case (3)}}}|\E[h_{c,k}(\bm{X}_{t_{1},k},\dots,\bm{X}_{t_{c},k})h_{c,k}(\bm{X}_{t_{c+1},k},\dots,\bm{X}_{t_{2c},k})]|\notag\\
&\hspace{4cm}=O(n^{2(c-1)})=O(n^{2c-1-\theta})
}
Finally, combining equations (\ref{case1m}), (\ref{case2m}) and (\ref{case3m}) yields the result. \hfill $\Box$


\subsubsection{Proof of (\ref{triangarray})} \label{sec:prooftriangarray}
We have by (\ref{dn2}) that
\al{
\hat f_{n,\xi}(\omega)
={}&\frac{1}{2\pi}\kabs\wk\Big\{\xi_k+\frac{m}{n-|k|}\sum_{t\in\mathcal{T}_k}h_{1,k}^{\xi}(\bm{X}_{t,k})\Big\}e^{-ik\omega}+O_{\Prob}\Big(r_nn^{-1/2-\theta/2}\Big).
}
Next, 
\al{
\hat f_{n,\tau}(\omega)-&\tilde f_{n,\tau}(\omega)
={}\frac{1}{2\pi}\kabs\wk\Big(\frac{m}{n-|k|}-\frac{m}{n}\Big)\sum_{t\in\mathcal{T}_k}h_{1,k}^{\xi}(\bm{X}_{t,k})e^{-ik\omega}\\&+\frac{1}{2\pi}\kabs\wk\frac{m}{n}\Big(\sum_{t\in\mathcal{T}_k}h_{1,k}^{\xi}(\bm{X}_{t,k})-\sum_{t=1}^{n}h_{1,k}^{\xi}(\bm{X}_{t,k})\Big)e^{-ik\omega}+O_{\Prob}\Big(r_nn^{-1/2-\theta/2}\Big)\\
={}&A_n+B_n+O_{\Prob}\Big(r_nn^{-1/2-\theta/2}\Big)
}
Similar arguments as in the proof of (\ref{dnlin}) yield
\al{
\E|A_n|\leq \frac{m}{2\pi}(2r_n+1)\Big|\frac{1}{n-r_n}-\frac{1}{n}\Big|\sup_{|k|\leq r_n}\E\Big|\sum_{t\in\mathcal{T}_k}h_{1,k}^{\xi}(\bm{X}_{t,k})\Big|=O\Big(\frac{r_n^2}{n^{3/2}}\Big),
}
where we have used that $|\frac{1}{n-r_n}-\frac{1}{n}|=O\Big(\frac{r_n}{n^2}\Big)$ and $\sup_{|k|\leq r_n}\E|\sum_{t\in\mathcal{T}_k}h_{1,k}^{\xi}(\bm{X}_{t,k})|=O(n^{1/2})$. Next,
\al{
B_n={}&-\frac{m}{2\pi}\sum_{\kip}\wk \frac{1}{n}\sum_{t=n-k+1}^nh_{1,k}^{\xi}(\bm{X}_{t,k})e^{-ik\omega}\\
&-\frac{m}{2\pi}\sum_{\kin}\wk \frac{1}{n}\sum_{t=1}^{|k|}h_{1,k}^{\xi}(\bm{X}_{t,k})e^{-ik\omega}\\
=:{}&B_{1,n}+B_{2,n}.
}
Note that by the stationarity of the process $\proc$,
\al{
B_{1,n}\eid-\frac{m}{2\pi}\sum_{\kip}\wk \frac{1}{n}\sum_{t=1}^kh_{1,k}^{\xi}(\bm{X}_{t,k})e^{-ik\omega}.
} 
Similarly as for $A_n$ we obtain
\al{
\E|B_{1,n}|\leq \frac{m}{2\pi}(2r_n+1)\frac{1}{n}\sup_{\kip}\E\Big|\sum_{t=1}^{k}h_{1,k}^{\xi}(\bm{X}_{t,k})\Big|=O\Big(\frac{r_n^{3/2}}{n}\Big)
} 
and analogously,
\al{
\E|B_{2,n}|\leq \frac{m}{2\pi}(2r_n+1)\frac{1}{n}\sup_{\kin}\E\Big|\sum_{t=1}^{|k|}h_{1,k}^{\xi}(\bm{X}_{t,k})\Big|=O\Big(\frac{r_n^{3/2}}{n}\Big).
}
Altogether, 
\al{
\hat f_{n,\tau}(\omega)-\tilde f_{n,\tau}(\omega)=O_{\Prob}\Big(r_nn^{-1/2-\theta/2}+r_n^{3/2}n^{-1}\Big).
}
This concludes the proof of~\eqref{triangarray}.{} \hfill $\Box$

\bigskip

\end{document}